\documentclass[11pt]{amsart}
\usepackage{amsmath,amssymb}
\usepackage{color}
\usepackage[all,curve,frame]{xy}

\usepackage[pdftex]{graphicx}
\usepackage{xcolor}

\newtheorem{defi}{Definition}[section]
\newtheorem{sinobservacion}[defi]{}
\newenvironment{sinob}{\begin{sinobservacion} \rm}{\end{sinobservacion} }
\newtheorem{coro}[defi]{Corollary}
\newtheorem{lem}[defi]{Lemma}
\newtheorem{rem}[defi]{Remark}
\newtheorem{prop}[defi]{Proposition}
\newtheorem{teo}[defi]{Theorem}
\newtheorem{ej}[defi]{Example}

\newcommand{\benu}{\begin{enumerate}}
\newcommand{\enu}{\end{enumerate}}

\begin{document}
\title[Enlargements of complexes of fixed size]
{Enlargements  of complexes of fixed size}

\author[C. Chaio]{Claudia Chaio}
\address[Claudia Chaio]{Centro Marplatense de Investigaciones Matem\'aticas, Facultad de Ciencias Exactas y
Naturales, Funes 3350, Universidad Nacional de Mar del Plata and CONICET,  Mar del
Plata, 7600, Argentina}
\email{cchaio@mdp.edu.ar}

\author[A. Gonz\'alez Chaio]{Alfredo Gonz\'alez Chaio}
\address[Alfredo Gonz\'alez Chaio]{Centro Marplatense de Investigaciones Matem\'aticas, Facultad de Ciencias Exactas y
Naturales, Funes 3350, Universidad Nacional de Mar del Plata,  Mar del
Plata, 7600, Argentina}
\email{agonzalezchaio@gmail.com}

\author[P. Suarez]{Pamela Suarez}
\address[Pamela Suarez]{Centro Marplatense de Investigaciones Matem\'aticas, Facultad de Ciencias Exactas y
Naturales, Funes 3350, Universidad Nacional de Mar del Plata,  Mar del
Plata, 7600, Argentina}
\email{psuarez@mdp.edu.ar}

\date{\today}


\maketitle

Corresponding Author: Claudia Chaio.

N/A (Data Availability Statement not applicable).

\begin{abstract} Let $A$ be an artin algebra. The aim of this work is to describe the enlargements of an indecomposable complex in $\mathbf{C}_{n}(\mbox{proj} \,A)$,  and to study the irreducible morphisms between them. Precisely,
we prove that any indecomposable complex in $\mathbf{C}_{[0,n]}(\mbox{proj} \,A)$  or in $\mathbf{C}_{n+1}(\mbox{proj} \,A)$ for $n$ a positive integer is a shift or an enlargement of an indecomposable complex in $\mathbf{C}_{n}(\mbox{proj} \,A)$. 
We also describe the entrances of the irreducible morphisms in $\mathbf{C}_{[0,n]}(\mbox{proj} \,A)$ between enlargements of an indecomposable complex $X$ in $\mathbf{C}_{n}(\mbox{proj} \,A)$.
\end{abstract}

\section*{Introduction}

Throughout this work, we consider $A$ an artin algebra, $\mbox{mod}\,A$ and $\mbox{proj}\,A$ the category of finitely generated $A$-modules and projective $A$-modules, respectively. 

In \cite{BSZ}, R. Bautista, M. J. Souto Salorio and  R. Zuazua defined the category of complexes of fixed size,
$\mathbf{C_n}(\mbox{proj}\, A)$  for $n \geq 2$, as the full subcategory of $\mathbf{C^{b}}({\rm mod}\,A)$ whose objects are the complexes
$X=\{(X^{i},d^{i})\}_{i \in \mathbb{Z}}$
such that $X^i \in  {\rm proj}\,A$ and $X^{i}=0$ if $i \notin \{1, \dots, n\}$.
Precisely, the above-mentioned authors defined and studied these categories in order to obtain information on the bounded derived category.

For each $n \geq 2$, the categories of complexes of fixed size are Krull-Schmidt categories, exact and with enough projective and injective objects, but they are not abelian categories.

{It is of interest to obtain information about objects and irreducible morphisms in $\mathbf{C_n}(\mbox{proj}\, A)$ whenever $n$ increases. In particular, the description of the category  $\mathbf{C_{[0,n]}}(\mbox{proj} \,A)$ in terms of the category  $\mathbf{C_{n}}(\mbox{proj} \,A)$  allows us to transfer information from the categories of complexes of fixed size to the bounded derived category.}

Following \cite{CGP6}, we say that an indecomposable complex $Z$ in $\mathbf{C_{[0,n]}}(\mbox{proj} \,A)$ is a left enlargement of an indecomposable complex $X$ in $\mathbf{C_n}(\mbox{proj}\, A)$ if $Z^{0} \neq 0$ and $X$ is a direct summand of $F(Z)$ in $\mathbf{C_{n}}(\mbox{proj} \,A)$.  Moreover, the above authors showed that enlargements are useful in determining whether an irreducible morphism in $\mathbf{C_{n}}(\mbox{proj} \,A)$ is still irreducible considered as a morphism in $\mathbf{C_{[0,n]}}(\mbox{proj} \,A)$.

In this work, we prove that any indecomposable complex in $\mathbf{C_{[0,n]}}(\mbox{proj} \,A)$ can be viewed as an enlargement or a shift of an indecomposable complex in $\mathbf{C_{n}}(\mbox{proj} \,A)$, see Theorem \ref{indecomposable}. Hence a complete description of the enlargements shields a complete description of the indecomposable objects of $\mathbf{C_{[0,n]}}(\mbox{proj} \,A)$ in terms of the indecomposable complexes of $\mathbf{C_{n}}(\mbox{proj} \,A)$. Furthermore, this description can also be used to describe irreducible morphisms. 

These results are expected to be helpful in finding  necessary and sufficient conditions over an irreducible morphism $f:X \rightarrow Y$  in $\mathbf{C_n}(\mbox{proj}\, A)$ to factor through an enlargement in $\mathbf{C_{[0,n]}}(\mbox{proj} \,A)$.

%{\color{blue} One of our aims is to find necessary and sufficient conditions over an irreducible morphism $f:X \rightarrow Y$  in $\mathbf{C_n}(\mbox{proj}\, A)$ to factor through an enlargement in $\mathbf{C_{[0,n]}}(\mbox{proj} \,A)$. For that reason, it shall be necessary to study the shape of the enlargements.} 

%\noindent{\bf Theorem A ([Theorem \ref{indecomposable}]).}  {\it Any indecomposable complex in $\mathbf{C_{[0,n]}}(\emph{proj} \,A)$ or in $\mathbf{C_{n+1}}(\emph{proj} \,A)$ is a shift  or an enlargement of an indecomposable complex in $\mathbf{C_{n}}(\emph{proj} \,A)$.}
%\vspace{.1in}
%Therefore, that is one of our reasons to develop the description of the enlargements of a given complex. 

{The main focus of this article is to study left enlargements (dually, right enlargements) and irreducible morphisms between them.}

In order to simplify the shape of an enlargement, we define the notion of diagonal complexes. Diagonal complexes are complexes $Z$ in $\mathbf{C_{[0,n]}}(\mbox{proj} \,A)$ where the differentials $d_{Z}^{i}$ for $i=1, \dots, n$ are diagonal matrices, see Definition \ref{diagonal}. Precisely, we prove Theorem A.
\vspace{.1in}

\noindent{\bf Theorem A ([Theorem \ref{diagonal}]).}  {\it Any non-injective complex in $\mathbf{C_{[0,n]}}(\emph{proj} \,A)$ is isomorphic to a
diagonal complex in $\mathbf{C_{[0,n]}}(\emph{proj} \,A)$.}
\vspace{.1in}

We also give conditions for the existence of a left enlargement, see Proposition \ref{existence}.

Given an indecomposable complex $X \in \mathbf{C}_{n}(\mbox{proj} \,A)$, we study the first entry of any left enlargement that appears in the almost split sequence in $\mathbf{C}_{[0,n]}(\mbox{proj} \,A)$  starting at $X$. This result is stated in Theorem B.
\vspace{.1in}

\noindent{\bf Theorem B ([Proposition \ref{prop-summads}]).}  {\it Let $X$ be an indecomposable complex in $\mathbf{C_{n}}(\emph{proj} \,A)$ and $f:X \rightarrow \overleftarrow{X}$ be an irreducible morphism in $\mathbf{C_{[0,n]}}(\emph{proj} \,A)$ where $\overleftarrow{X}$ is a left enlargement of $X$. Assume that $\overleftarrow{X}$ is a complex as follows $${\xymatrix {
Z^{0}  \ar[rr]^{\left(\begin{array}{cc} \widetilde{d}_{X}^{0} \\  \widetilde{d}_{Z}^{0}  \end{array} \right)}&&  X^{1}\oplus Z^{1}
 \ar[rr]^{\left(\begin{array}{cc} d_{X}^{1} &0 \\  0 & d_{Z}^{1}  \end{array} \right)}&& X^{2}\oplus Z^{2}   \ar[rr]^{\left(\begin{array}{cc} d_{X}^{2} &0 \\ 0 & d_{Z}^{2}  \end{array} \right)}&& \dots  \ar[rr]^{\left(\begin{array}{cc} d_{X}^{n-1} &0 \\  0 & d_{Z}^{n-1}  \end{array} \right)} &&  X^{n}\oplus Z^{n}}}.$$ Then $Z^{0}$ is a direct summand of the projective cover of $\emph{soc}(\emph{Ker}\, {d}_{X}^{1})$.}
\vspace{.1in}

In the last section of this work, we describe the entrances of the irreducible morphisms in $\mathbf{C_{[0,n]}}(\mbox{proj} \,A)$ between left enlargements of an indecomposable complex $X$ in $\mathbf{C}_{n}(\mbox{proj} \,A)$, see Theorem \ref{enlarg-enlarg}. { Moreover, we give a method to decide the type of irreducible morphisms in the sense of H. Giraldo and H. Merklen.
\vspace{.2in}

The paper is organized as follows. In Section 1, we introduce some necessary preliminaries. In Section 2, we give a description of indecomposable complexes and irreducible morphisms of $\mathbf{C_{[0,n]}}(\mbox{proj} \,A)$ in terms of indecomposable complexes and irreducible morphisms of ${\mathbf {C_n}}({\rm proj}\,A)$. We also prove Theorem A. Section 3, is dedicated to  study  the existence of particular enlargements. Furthermore, we give necessary and sufficient conditions to determine the indecomposability of some special complexes. Section 4, is dedicated to prove Theorem B. Finally, in Section 5, we study irreducible morphisms between enlargements of the same indecomposable complex and, we give a method to decide the type of such morphisms.}

\vspace{.2in}

\thanks{The authors thankfully acknowledge partial support from Universidad Nacional de Mar del Plata, Argentina.
The first author is a researcher of CONICET Argentina.}

\section{Preliminaries}

We start given some notations and basic facts useful for the subsequent sections.

\begin{sinob}

Let $A$ be an~artin algebra, ${\rm proj}\,A$ and  ${\rm inj}\,A$  be the full subcategories of ${\rm mod}\,A$ whose objects are the projective and the injective $A$-modules, respectively.  $\mbox{add}\,X$ means the full subcategory of ${\rm mod}\,A$ whose objects are the direct sums of direct summands of $X$.

Let $\mathbb{Z}$ be the set of integer numbers. For $m,n \in \mathbb{Z}$, with $n-m \geq 1$,  consider $[m,n]$ an interval in
$\mathbb{Z}$,  $\mathbf{C_{[m,n]}}({\rm proj}\, A)$  is the full subcategory of $\mathbf{C^{b}}({\rm mod}\,A)$ whose objects are the
complexes $X=\{(X^{i},d^{i})\}_{i \in \mathbb{Z}}$  such that $X^i \in  {\rm proj}\,A$ and $X^{i}=0$ if $i\notin \{m, \dots, n\}$.
We denote $\mathbf{C_{[1,n]}}({\rm proj}\, A)$ by $\mathbf{C_{n}}(\mbox{proj} \,A)$. %Similarly, we can define  the category $\mathbf{C_{[m,n]}}({\rm inj}\, A)$.

%Let $m,n,m',n'\in \mathbb{Z}$  and $[m,n] \subset [m',n']$. We may consider $\mathbf{C_{[m,n]}}({\rm proj}\, A)  \subset \mathbf{C_{[m',n']}}({\rm proj}\, A)$ and moreover, if $n-m=n'-m'$ then the categories $\mathbf{C_{[m,n]}}({\rm proj}\, A)$  and $\mathbf{C_{[m',n']}}({\rm proj}\, A)$ are isomorphic. Similarly, the categories $\mathbf{C_{[m,n]}}({\rm inj}\, A)$  and $\mathbf{C_{[m',n']}}({\rm inj}\, A)$ are also isomorphic.

Consider $\mathbf{C^{\leq n}}({\rm proj}\,A)$ the full subcategory of $\mathbf{C}({\rm proj}\,A)$ whose objects are those $X \in \mathbf{C}({\rm proj}\,A)$
such that $X^i=0$ for $i>n$. Similarly, $\mathbf{C^{\geq n}}({\rm proj}\,A)$ is the full subcategory of $\mathbf{C}({\rm proj}\,A)$ whose objects are
those $X \in \mathbf{C}({\rm proj}\,A)$ such that $X^i=0$ for $i < n$.

Let $F: \mathbf{C^{\leq n}}({\rm proj}\,A)\rightarrow \mathbf{C_n}({\rm proj} \,A)$ and
$G: \mathbf{C^{\geq 1}}({\rm proj}\,A) \rightarrow \mathbf{C_n}({\rm proj} \,A)$ be the left and right hard truncation functors, respectively.
\end{sinob}

\begin{sinob}\label{K.S}
Let $R$ be a fixed commutative artinian ring, and $\mathcal{A}$ be a Krull-Schmidt $R$-category,
that is a Hom-finite additive R-category in which idempotents split.

The Jacobson radical $\Re(\mathcal{A})$ of $\mathcal{A}$ is the ideal generated by the non-isomorphisms between indecomposable
objects. For $n > 1$, $\Re^{n}(\mathcal{A})$ denotes the ideal generated by the morphisms
which are composition of $n$ morphisms lying in the radical and $\Re^{\infty}(\mathcal{A}) =
\cap_{n \geq 1}\Re^{n}(\mathcal{A})$ the infinite radical of $\mathcal{A}$.
\vspace{.05in}

Let $f : X \rightarrow Y$ be a morphism in $\mathcal{A}$. We recall that $f$ is irreducible
if $f$ is neither a section nor a retraction and whenever $f = h g$ then either $g$ is a section or $h$ is a retraction.

If $X, Y$ are indecomposable objects in $\mathcal{A}$, then $$\mbox{Irr}(X, Y) = \Re(X, Y)/\Re^2(X, Y)$$ is a
$k_{X}-k_{Y}$-bimodule, where $k_Z = \mbox{End}\,(Z)/ \Re(\mbox{End}\,(Z))$.
\end{sinob}

\begin{sinob}\label{GM}
By  \cite[Corollary 2, Proposition 3]{GM:08} if
  $f = \{f^i\}_{i=1}^n:X \rightarrow Y$
is an irreducible
morphism in $\mathbf{C_n}({\rm proj} \,A)$ then, one of the
following statements hold.
\begin{enumerate}
\item{(sec).} For each $i \in \{1, \dots, n\}$, the morphisms $f^i$ are sections in ${\rm proj} \,A.$
\item{(ret).} For each $i \in \{1, \dots, n\}$, the morphisms $f^i$ are retractions in ${\rm proj} \,A.$
\item{(ret-irr-sec).} There exists $i \in \{1, \dots, n\}$ such that $f^i$ is irreducible in ${\rm proj}\,A$, the morphisms $f^j$ are sections for
all $j>i$ and the morphisms $f^j$ are retractions for $j<i.$
\end{enumerate}
\end{sinob}

{
\begin{sinob}
We transcribe Proposition 2.2 stated in \cite{BS}, which shall be very useful for subsequent sections.

\begin{prop} \label{sec-ret-irr} \cite[Proposition 2.2]{BS}. Let $k$ be a commutative artinian ring and $\mathcal{A}$ be an additive $k$-category. 
Let $f:X \rightarrow Y$ be an irreducible morphism in $C_{J}(\mathcal{A})$ with $J$ an interval or $J=  \mathbb{Z}$ and $I$ an interval contained in $J$. Then, $f_{I}:X_{I} \rightarrow Y_{I}$ is a section, a retraction or an irreducible morphism in $C_{I}(\mathcal{A})$.
\end{prop}
\end{sinob}}

\begin{sinob}\label{resultados-cono}

The mapping cone of a morphism of complexes  $f : X \rightarrow Y$ is the complex  $C(f)= \{((X[1])^{i} \oplus Y^{i}, d_{C(f)}^{i})\}_{i \in \mathbb{Z}}$
with differential $d_{C(f)}^{i}= \left(
                    \begin{array}{cc}
                      -d_X^{i+1} & 0 \\
                      f^{i+1}& d_Y^{i} \\
                    \end{array}
                  \right)$.
\end{sinob}

\begin{sinob}\label{extended}
We recall the definition given in \cite{CGP2} of   when a complex  in $\mathbf{C_{n}}(\mbox{proj} \,A)$ can be extended  to the right or to the left.

\begin{defi} \label{exten} Let $X=\{(X^i, d^i)\}_{i\in \mathbb{Z}}$ be an indecomposable complex in $\mathbf{C_{n}}(\emph{proj} \,A)$.
We say that $X$ can be extended to the right if there  is $Q^{n+1} \in \emph{proj} \,A$ and a non-zero morphism $d'^{n}: X^{n} \rightarrow Q^{n+1}$
such that $d'^{n} d^{n-1}=0$.

Dually, we say that $X$ can be extended to the left if there is  $Q^{0} \in \emph{proj} \,A$ and a non-zero morphism $d'^{0}: Q^{0} \rightarrow X^{1}$
such that $d^{1} d'^{0}=0$.
\end{defi}

\end{sinob}

\begin{sinob}

% Inspired in Lemma \ref{extension},  we recall a more general definition than the one given in Definition \ref{exten} of when a complex
% in $\mathbf{C_{n}}(\mbox{proj} \,A)$ can be extended to right or to the left. Precisely, we shall not ask the complex $X$ to be indecomposable.

% \begin{defi} \label{exten-nueva 2}  Let $X=\{(X^i, d^i)\}_{i\in \mathbb{Z}}$ be a complex in $\mathbf{C_{n}}(\emph{proj} \,A)$. We say that $X$
% can be extended to the left if $d^1$ is not a monomorphism. Dually, we say that
% $X$ can be extended to the right if $\emph{Hom}_A( \emph{coker}\, d^{n-1}, A) \neq 0$.
% \end{defi}

Below, we recall the definition of a left and a right enlargement of an indecomposable complex  given in \cite{CGP6}.

\begin{defi} \label{enlargment-left} Let $X$ be an indecomposable complex in $\mathbf{C_{n}}(\emph{proj} \,A)$. We say that an
indecomposable complex $Z$ in $\mathbf{C_{[0,n]}}(\emph{proj} \,A)$ is a left enlargement of $X$  if $Z^{0} \neq 0$ and $X$ is a
direct summand of $F(Z)$ in  $\mathbf{C_{n}}(\emph{proj} \,A)$. We denote a left enlargement of $X$  by $\overleftarrow{X}$.
 \end{defi}

\begin{defi} \label{enlargment-right} Let $X$ be an indecomposable complex in $\mathbf{C_{n}}(\emph{proj} \,A)$. We say that an indecomposable
complex $Z$ in $\mathbf{C_{n+1}}(\emph{proj} \,A)$ is a right enlargement of $X$ if $Z^{n+1} \neq 0$ and  $X$ is a direct summand of $G(Z)$ in $\mathbf{C_{n}}(\emph{proj} \,A)$.  We denote a  right enlargement of $X$  by $\overrightarrow{X}$.
 \end{defi}

Let $X$ be an indecomposable complex in $\mathbf{C_{n}}(\mbox{proj} \,A)$. We say that  a left enlargement $Z$ of $X$ is a left extension of $X$ if $F(Z)=X$ in  $\mathbf{{C}_{n}}(\mbox{proj} \,A)$.

Dually, we say that a right enlargement $Z$ of $X$ in $\mathbf{C_{n+1}}(\mbox{proj} \,A)$ is a right extension of $X$  if
$G(Z)=X$ in  $\mathbf{C_{n}}(\mbox{proj} \,A)$.
\vspace{.05in}

The next lemma shall be useful for our considerations.

\begin{lem} \label{lema-previo}\cite[Lemma 2.4]{CGP6} Let $X$ be an indecomposable complex  in $\mathbf{C_{n}}(\emph{proj} \,A)$. The following conditions hold.
\begin{enumerate}
\item[(a)] The direct summands of the codomain of the left minimal almost split morphism starting in $X$ in $\mathbf{C_{[0,n]}}(\emph{proj} \,A)$
are of the form $\overleftarrow{X}$ or $W_i$ where for all $i$, $W_i$ is a direct summand of the codomain of the left minimal almost split morphism starting in $X$ in $\mathbf{C_{n}}(\emph{proj} \,A)$.
\item[(b)] The direct summands of the domain of the right minimal almost split morphism ending in $X$ in $\mathbf{C_{n+1}}(\emph{proj} \,A)$ are of the
form $\overrightarrow{X}$ or $W_i$ where for all $i$, $W_i$ is a direct summand of the domain of the  right minimal almost split morphism ending in $X$
in $\mathbf{C_{n}}(\emph{proj} \,A)$.
\end{enumerate}
\end{lem}
\end{sinob}

\begin{sinob} The next definition and results from \cite{BS} shall be fundamental to analyze how are the entrances of irreducible
morphisms in $\mathbf{C_{n+1}}(\mbox{proj} \,A)$ between enlargements of an indecomposable complex $X \in \mathbf{C_{n}}(\mbox{proj} \,A)$.

\begin{defi} A morphism $g: Z \rightarrow W$ in a Krull-Schmidt category $\mathcal{A}$  is called radical if for any section $\sigma: Z_1 \rightarrow Z$ and
any retraction $\tau: W \rightarrow W_1$, the composition $\tau g \sigma$ is not an isomorphism.
\end{defi}

Observe that a morphism in a  Krull-Schmidt category $\mathcal{A}$  is radical if and only if it belongs to the radical of $\mathcal{A}$.

\begin{prop} \cite[Proposition 2.16]{BS} \label{BS1} Let $\mathcal{A}$ be a Krull-Schmidt category, then any morphism $f:X \rightarrow  Y$ is isomorphic to a morphism of the form:

$$f=\left(
      \begin{array}{cc}
        f_1 & 0 \\
        0 & f_2 \\
      \end{array}
    \right) : X_1 \oplus X_2 \rightarrow  Y_1 \oplus Y_2$$

\noindent with $f_1$ an isomorphism and $f_2$ a radical morphism.
\end{prop}

\begin{prop} \cite[Proposition 2.17]{BS} \label{BS2} Let $\mathcal{A}$ be a Krull-Schmidt category. A morphism

$$f=\left(
      \begin{array}{cc}
        f_1 & 0 \\
        0 & f_2 \\
      \end{array}
    \right) : X_1 \oplus X_2 \rightarrow  Y_1 \oplus Y_2,$$

\noindent with $f_1$ an isomorphism and $f_2$ a radical morphism, is irreducible if and only if $f_2$ is irreducible.
\end{prop}

\begin{prop} \cite[Proposition 2.18]{BS} \label{BS3} Let $f:X \rightarrow  Y$  be an irreducible morphism in the radical of a  Krull-Schmidt
category $\mathcal{A}$. Then if $p:Y \rightarrow  Z$ is a non-zero retraction, the morphism $pf: Y \rightarrow  Z$ is irreducible.
Similarly, if $s:W\rightarrow  X$ is a non-zero section then $fs: W \rightarrow  Y$ is an irreducible morphism.
\end{prop}
\end{sinob}

Notation: We say that $f:X\rightarrow Y$ is a radical-irreducible morphism if it is a morphism in the radical of a Krull-Schmidt category which is also irreducible. 

\section{Diagonal complexes}

{In this section, we define the notion of  diagonal complexes. We show that any non-injective indecomposable complex in $\mathbf{C_{[0,n]}}(\mbox{proj} \,A)$ is isomorphic to a diagonal complex in
$\mathbf{C_{[0,n]}}(\mbox{proj} \,A)$.

We also prove that indecomposable complexes and irreducible morphisms in $\mathbf{C_{[0,n]}}(\mbox{proj} \,A)$ can be described in terms of indecomposable complexes and irreducible morphisms in $\mathbf{C_{n}}(\mbox{proj} \,A)$.}
\vspace{.1in}
 
We start with the following definition.

\begin{defi} \label{diagonal} Let $Z \in \mathbf{C_{[0,n]}}(\emph{proj} \,A)$ be the following complex:

$$Z^{0}\stackrel{\left(\begin{array}{cc} d_{X}^{0} \\  d_{Z}^{0}  \end{array} \right)} \longrightarrow X^{1} \oplus Y^{1}
\stackrel{\left(\begin{array}{cc} d_{X}^{1} & a_1\\  b_1 & d_{Y}^{1}  \end{array} \right)}
 \longrightarrow X^{2} \oplus Y^{2}\stackrel{\left(\begin{array}{cc} d_{X}^{2} & a_2\\  b_2 & d_{Y}^{2}  \end{array} \right)}  \longrightarrow
 \dots \stackrel{\left(\begin{array}{cc} d_{X}^{n-1} & a_{n-1}\\  b_{n-1} & d_{Y}^{n-1}  \end{array} \right)}  \longrightarrow X^{n} \oplus Y^{n}$$

\noindent with ${X:\xymatrix {  X^{1} \ar[r]^{d_{X}^{1}} & \dots \ar[r]^{d_{X}^{n-1}} & X^{n}}}$  an indecomposable complex  in $\mathbf{C_{n}}(\emph{proj} \,A)$.% and $Z^{0} \neq 0$. 
We say that $Z$ is a {\bf diagonal complex} if $a_i=0$ and $b_i=0$ for all  $i=1, \dots, n-1$.
\end{defi}

 \vspace{.1in}

Now, we are in a position to prove one of the main results of this section.

\begin{teo} \label{description} Any non-injective complex in $\mathbf{C_{[0,n]}}(\emph{proj} \,A)$ is isomorphic to a diagonal complex in
$\mathbf{C_{[0,n]}}(\emph{proj} \,A)$.
\end{teo}

%Let $X$ be an indecomposable complex in $\mathbf{C_{n}}(\emph{proj} \,A)$ and $\overleftarrow{X}: P_0 \rightarrow X^{1} \oplus Y^{1}\rightarrow X^{2} \oplus Y^{2}\rightarrow \dots $ be a right  enlargement of $X$ in $\mathbf{C_{[0,n]}}(\emph{proj} \,A)$ where $d^{i}= $. Assume that  there is an irreducible morphism from $X$ to $\overleftarrow{X}$ in $\mathbf{C_{[0,n]}}(\emph{proj} \,A)$. Then  for all $i$, $b_i=0$.

\begin{proof}  Let $Z$ be a complex in $\mathbf{C_{[0,n]}}(\mbox{proj} \,A)$ and ${X:\xymatrix {  X^{1} \ar[r]^{d_{X}^{1}} & \dots \ar[r]^{d_{X}^{n-1}} & X^{n}}}$ be an indecomposable direct summand of $F(Z)$. Then $$Z: Z^{0}\stackrel{\left(\begin{array}{cc} \widetilde{d}_{X}^{0} \\  \widetilde{d}_{Y}^{0}  \end{array} \right)} \longrightarrow X^{1} \oplus Y^{1} \stackrel{\left(\begin{array}{cc} d_{X}^{1} & a_1\\  b_1 & d_{Y}^{1}  \end{array} \right)}
 \longrightarrow X^{2} \oplus Y^{2}\stackrel{\left(\begin{array}{cc} d_{X}^{2} & a_2\\  b_2 & d_{Y}^{2}  \end{array} \right)}  \longrightarrow  \dots \stackrel{\left(\begin{array}{cc} d_{X}^{n-1} & a_{n-1}\\  b_{n-1} & d_{Y}^{n-1}  \end{array} \right)}  \longrightarrow X^{n} \oplus Y^{n}$$
 \vspace{.1in}

\noindent with $$F(Z): X^{1} \oplus Y^{1} \stackrel{\left(\begin{array}{cc} d_{X}^{1} & a_1\\  b_1 & d_{Y}^{1}  \end{array} \right)}
 \longrightarrow X^{2} \oplus Y^{2}\stackrel{\left(\begin{array}{cc} d_{X}^{2} & a_2\\  b_2 & d_{Y}^{2}  \end{array} \right)}  \longrightarrow  \dots \stackrel{\left(\begin{array}{cc} d_{X}^{n-1} & a_{n-1}\\  b_{n-1} & d_{Y}^{n-1}  \end{array} \right)}  \longrightarrow X^{n} \oplus Y^{n}.$$
 \vspace{.1in}

Since $X$ is a direct summand of $F(Z)$, then there is a section $h: X \rightarrow F(Z).$ Without loss of generality we may assume that each $h^{i}= {\left(\begin{array}{cc} 1 \\ 0 \end{array} \right)}$ for $i=1, \dots, n$. Therefore, we have the following  induced morphism of complexes:

\begin{equation}\label{diagrama}
\tiny{\xymatrix {  0 \ar[dd]^{0} \ar[rr]^{0}&&  X^{1}   \ar[dd]^{\left(\begin{array}{cc} 1 \\ 0 \end{array} \right)} \ar[rr]^{d_{X}^{1}}&& X^{2} \ar[dd]^{\left(\begin{array}{cc} 1 \\ 0 \end{array} \right)}\ar[rr]^{d_{X}^{2}} &&
\dots \ar[r]& X^{n-1} \ar[dd]^{\left(\begin{array}{cc} 1 \\ 0 \end{array} \right)}\ar[rr]^{d_{X}^{n-1}} && X^{n} \ar[dd]^{\left(\begin{array}{cc} 1 \\ 0 \end{array} \right)} &&\\ &&  &&   & & &&\\
Z ^{0} \ar[rr]^{\left(\begin{array}{cc} \widetilde{d}_{X}^{0} \\ \widetilde{d}_{Y}^{0}  \end{array} \right)} && X^{1}\oplus Y^{1} \ar[rr]^{\left(\begin{array}{cc} d_{X}^{1} & a_1\\  b_1 & d_{Y}^{1}  \end{array} \right)} && X^{2}\oplus Y^{2} \ar[rr]^{\left(\begin{array}{cc} d_{X}^{2} & a_2\\  b_2 & d_{Y}^{2}  \end{array} \right)} && \dots \ar[r] & X^{n-1}\oplus Y^{n-1}\ar[rr]^{\left(\begin{array}{cc} d_{X}^{n-1}& a_{n-1}\\  b_{n-1} & d_{Y}^{n-1}  \end{array} \right)} && X^{n}\oplus Y^{n}&&.}}
\end{equation}
\vspace{.1in}

Hence, for all $i=1, \dots, n-1$ we have that

$$ {\left(\begin{array}{cc} d_{X}^{i} & a_i\\  b_i & d_{Y}^{i}  \end{array} \right)}{\left(\begin{array}{cc} 1 \\ 0 \end{array} \right)}= {\left(\begin{array}{cc} 1 \\ 0 \end{array} \right)} d_{X}^{i}$$ getting that $ b_i =0$.

Now, since

$${\left(\begin{array}{cc} d_{X}^{i+1} & a_{i+1}\\  o & d_{Y}^{i+1}  \end{array} \right)} {\left(\begin{array}{cc} d_{X}^{i} & a_i\\ 0 & d_{Y}^{i}  \end{array} \right)}= 0$$ %{\left(\begin{array}{cc}0 & 0\\ 0& 0  \end{array} \right)}$$
\noindent because $F(Z)$ is a complex, then we get that

$$ d_{X}^{i+1} a_i= - a_{i+1} d_{Y}^{i}$$
%\vspace{.1in}

\noindent for all $i=1, \dots, n-1$.

On the other hand, since $h: X \rightarrow F(Z)$ is a section then there is a morphism ${\left(\begin{array}{cc} 1 & c_{i}\end{array} \right)} :X^{i}\oplus Y^{i} \rightarrow X^{i}$ such that ${\left(\begin{array}{cc} 1 & c_{i}\end{array} \right)}{\left(\begin{array}{cc} 1 \\ 0 \end{array} \right)}= 1_{X^{i}},$ for $i=1, \dots, n$. Moreover, ${\left(\begin{array}{cc} 1 & c\end{array} \right)}:F(Z) \rightarrow X$ is a morphism of complexes. Then

$${\left(\begin{array}{cc} 1 & c_{i+1}\end{array} \right)}
{\left(\begin{array}{cc} d_{X}^{i} & a_{i}\\  o & d_{Y}^{i}  \end{array} \right)} = {\left(\begin{array}{cc} d_{X}^{i} & d_{X}^{i} c_{i}  \end{array} \right)}, $$

\noindent getting that %$\widetilde{d}_{X}^{0}=-c_{1} \widetilde{d}_{Y}^{0}$ and 
$a_i= d_{X}^{i}c_{i}- c_{i+1} d_{Y}^{i}$,  for $i=1, \dots, n-1$.
%In particular, $a_1= d_{X}^{i}c_{1}- c_{2} d_{Y}^{1}$.

Now, consider $d^{\prime}=  \widetilde{d}_{X}^{0} + c_1  \widetilde{d}_{Y}^{0}$ and $\widetilde{Z}$ as follows:

 $$\widetilde{Z}: Z^{0}\stackrel{\left(\begin{array}{cc} d^{\prime} \\  \widetilde{d}_{Y}^{0}  \end{array} \right)} \longrightarrow X^{1} \oplus Y^{1} \stackrel{\left(\begin{array}{cc} d_{X}^{1} & 0\\ 0 & d_{Y}^{1}  \end{array} \right)}
 \longrightarrow X^{2} \oplus Y^{2}\stackrel{\left(\begin{array}{cc} d_{X}^{2} & 0\\  0 & d_{Y}^{2}  \end{array} \right)}  \longrightarrow  \dots \stackrel{\left(\begin{array}{cc} d_{X}^{n-1} & 0\\  0 & d_{Y}^{n-1}  \end{array} \right)}  \longrightarrow X^{n} \oplus Y^{n}.$$

 \vspace{.1in}

We affirm that $\widetilde{Z}$ is a complex. In fact, since $X$ and $Y$ are complexes then $d_{X}^{i+1} d_{X}^{i}=0$ and $ d_{Y}^{i+1} d_{Y}^{i}=0$ for $i=1, \dots, n-1$. It is left to prove that
${\left(\begin{array}{cc} d_{X}^{1} & 0\\  0 & d_{Y}^{1}  \end{array} \right)} {\left(\begin{array}{c} d^{\prime} \\ \widetilde{d}_{Y}^{0}  \end{array} \right)}=0.$

Since ${\left(\begin{array}{cc} d_{X}^{1} & a_1\\  0 & d_{Y}^{1}  \end{array} \right)} {\left(\begin{array}{c} \widetilde{d}_{X}^{0}\\ \widetilde{d}_{Y}^{0}  \end{array} \right)}=0$ then $d_{X}^{1}\widetilde{d}_{X}^{0}+ a_1 \widetilde{d}_{Y}^{0}=0$ and $ d_{Y}^{1} \widetilde{d}_{Y}^{0}=0$. Then it is only left to prove that $d_{X}^{1}d' =0$. In fact, $d_{X}^{1}d' = d_{X}^{1}(\widetilde{d}_{X}^{0} + c_1  \widetilde{d}_{Y}^{0})$. Furthermore, since $a_1= d_{X}^{i}c_{1}- c_{2} d_{Y}^{1}$ then $a_1\widetilde{d}_{Y}^{0}= d_{X}^{i}c_{1}\widetilde{d}_{Y}^{0}.$ Hence, $d_{X}^{1}d' =d_{X}^{1}\widetilde{d}_{X}^{0} + a_1  \widetilde{d}_{Y}^{0}=0$, proving that $\widetilde{Z}$ is a complex.

Finally, we prove that $\widetilde{Z}$ is isomorphic to $Z$. Consider the morphism $g=\{g^{i}\}_{i=0}^{n}$, given by the following diagram:

$${\tiny {\xymatrix { Z^{0} \ar[dd]^{1} \ar[rr]\ar[rr]^{\left(\begin{array}{cc} \widetilde{d}_{X}^{0} \\  \widetilde{d}_{Y}^{0}  \end{array} \right)}&&  X^{1}\oplus Y^{1}  \ar[dd]^{\left(\begin{array}{cc} 1 & c_1\\ 0 & id\end{array} \right)} \ar[rr]^{d_{X}^{1}}&& X^{2}\oplus Y^{2} \ar[dd]^{\left(\begin{array}{cc} 1 &c_2 \\ 0 &1 \end{array} \right)}\ar[rr] &&
\dots \ar[rr]&&  X^{n} \oplus Y^{n} \ar[dd]^{\left(\begin{array}{cc} 1 & c_n \\ 0 & 1\end{array} \right)} &&\\ &&  &&    &&\\
Z ^{0} \ar[rr]^{\left(\begin{array}{cc} d' \\  \widetilde{d}_{Y}^{0}  \end{array} \right)} && X^{1}\oplus Y^{1} \ar[rr]^{\left(\begin{array}{cc} d_{X}^{1} & a_1\\  b_1 & d_{Z}^{1}  \end{array} \right)} && X^{2}\oplus Y^{2} \ar[rr] && \dots \ar[rr] && X^{n}\oplus Y^{n}&&.}}}$$
\vspace{.1in}
\noindent where $g^{i}={\left(\begin{array}{cc} 1 & c_i \\ 0 & 1\end{array} \right)}$ for $i= 1, \dots, n$ and $g^{0}=id$.

Clearly, all the morphisms $g^{i}$ for $i=0, \dots, n$ are isomorphisms. It is not hard to see that each square commutes. Therefore, we conclude that $\widetilde{Z}$ is isomorphic to $Z$, proving that each complex is isomorphic to a diagonal complex.
\end{proof}

\begin{rem}
\emph{
Consider the following complex: 
$$Z: Z^0 \to Z^1 \to  \dots \to Z^n.$$ Let $F(Z)=\oplus^m_{i=1} X_i$ with every $X_i$ indecomposable.  Applying Theorem \ref{description} to first direct summand $X_1$ we get that $Z$ is isomorphic to:
$$Z:  Z^{0}\stackrel{\left(\begin{array}{cc} \widetilde{d}_{X_1}^{0} \\  \widetilde{d}_{Y}^{0}  \end{array} \right)} \longrightarrow X_1^{1} \oplus Y^{1} \stackrel{\left(\begin{array}{cc} d_{X_1}^{1} & 0 \\  0 & d_{Y}^{1}  \end{array} \right)}  \longrightarrow X_1^{2} \oplus Y^{2}\stackrel{\left(\begin{array}{cc} d_{X_1}^{2} & 0\\  0& d_{Y}^{2}  \end{array} \right)}  \longrightarrow  \dots \stackrel{\left(\begin{array}{cc} d_{X_1}^{n-1} & 0\\  0 & d_{Y}^{n-1}  \end{array} \right)}  \longrightarrow X_1^{n} \oplus Y^{n}.$$
 Successively applying Theorem \ref{description} in each indecomposable summand $X_i:   X_i^{1}  \stackrel{\widetilde{d}_{X_i}^{1}} \longrightarrow X_i^{2} \stackrel{\widetilde{d}_{X_i}^{2} }  \longrightarrow  \dots \stackrel{\widetilde{d}_{X_i}^{n-1} } \longrightarrow X_i^{n}$ of  $Y=\oplus^{m}_{i=2} X_i$, we obtain a complex $\widetilde{Z}$ which is isomorphic to $Z$. Precisely,  $\widetilde{Z}$ is as follow}
\emph{ 
 $$\widetilde{Z}:  Z^{0}\stackrel{\widetilde{d}_{Z}^{0} } \longrightarrow X_1^{1} \oplus X_2^{1} \dots X_m^{1} \stackrel{\widetilde{d}_{Z}^{1} }
 \longrightarrow X_1^{2} \oplus X_2^{2} \dots X_m^{2}\stackrel{\widetilde{d}_{Z}^{2} }  \longrightarrow  \dots \stackrel{\widetilde{d}_{Z}^{n-1} } \longrightarrow X_1^{n} \oplus X_2^{n} \dots X_m^{n}$$}

\emph{\noindent where the differentials in each entry are diagonal matrices as we described below:}

\emph{ $$\widetilde{d}_{Z}^{i}=\left(\begin{array}{cccc} d_{X_1}^{i} & 0 & \dots & 0\\  0 & d_{X_2}^{i}  &\dots & 0 \\ &   &\ddots &\\ 0 & 0  &\dots & d_{X_m}^{i}  \end{array} \right).$$}
 
\emph{We denote such a complex with the following expression ${\xymatrix {\widetilde{Z}:  Z^{0}  \ar[r] & X_1 \oplus \dots \oplus X_{m}  }}$.}
\end{rem}

The aim of the rest of this section is to show that indecomposable complexes and irreducible morphisms in $\mathbf{C_{[0,n]}}(\mbox{proj} \,A)$ can be described in terms of indecomposable complexes and irreducible morphisms in $\mathbf{C_{n}}(\mbox{proj} \,A)$.
\vspace{.1in}

We start with the following result which describes the indecomposable complexes of $\mathbf{C_{[0,n]}}(\mbox{proj} \,A)$.

\begin{teo}\label{indecomposable} 
Any indecomposable complex in $\mathbf{C_{[0,n]}}(\emph{proj} \,A)$ or in $\mathbf{C_{n+1}}(\emph{proj} \,A)$
is a shift or an enlargement of an indecomposable complex in $\mathbf{C_{n}}(\emph{proj} \,A)$.
\end{teo}

\begin{proof}  Consider $X$ an indecomposable complex in $\mathbf{C_{[0,n]}}(\mbox{proj} \,A)$. If $X^{0}=0$ then $F(X)=X$ where $X$ is an indecomposable
complex in $\mathbf{C_{n}}(\mbox{proj} \,A)$. Then $X=X[0]$ is a shift  of  $F(X)$. If $X^{n}=0$, then $X[-1] \in \mathbf{C_{n}}(\mbox{proj} \,A)$ and it
is an indecomposable complex in $\mathbf{C_{n}}(\mbox{proj} \,A)$. Hence $X= X[-1][1]$ is the shift of an indecomposable complex in $\mathbf{C_{n}}(\mbox{proj} \,A)$.

Now, assume that $X$ is the complex
${\xymatrix { X^{0}   \ar[r]^{d_{X}^{0}}& X^{1} \ar[r]^{d_{X}^{1}} & \dots \ar[r]^{d_{X}^{n-1}} & X^{n} }}\in \mathbf{C_{[0,n]}}(\mbox{proj} \,A)$
with $X^{0}\neq 0$ and $X^{n}\neq 0$. Then all entries $X^{i}$  of the complex $X$ are non-zero because $X$ is indecomposable. Furthermore, $F(X)$ is a
complex  as follows
$${\xymatrix {F(X):  X^{1}   \ar[r]^{d_{X}^{1}} & X^{2} \ar[r]^{d_{X}^{2}} & \dots \ar[r]^{d_{X}^{n-1}} & X^{n}}}\in \mathbf{C_{n}}(\mbox{proj} \,A).$$
%\noindent where each $X^{i}= \oplus_{j=1}^{n_j} Y^{i}_j$ with $Y^{i}_j$ indecomposable complexes in $\mathbf{C_{n}}(\mbox{proj} \,A)$,
Therefore, $F(X)=\oplus_{i=1}^{n} Y_i$ with $Y_i$ indecomposable complexes in $\mathbf{C_{n}}(\mbox{proj} \,A)$,
proving that $X= \overleftarrow{Y_i}$ for any of the direct summands of $F(X)$.

A similar argument gives us the result for a complex in $\mathbf{C_{n+1}}(\mbox{proj} \,A)$, which proves the statement.
\end{proof}

Finally, we describe the irreducible morphisms in
$\mathbf{C_{[0,n]}}(\mbox{proj} \,A)$ in terms of the irreducible of $\mathbf{C_{n}}(\mbox{proj} \,A)$.  
\vspace{.1in}

First, we prove a technical lemma useful to give a description of the irreducible morphisms between indecomposable complexes in terms of enlargements.  

\begin{lem}\label{noirred} Let $Y$ be an indecomposable complex in $\mathbf{C_{n}}(\emph{proj} \,A)$. There are not irreducible morphisms from any left enlargement of $Y$ to $Y$ in $\mathbf{C_{[0,n]}}(\emph{proj} \,A)$. 
\end{lem}
\begin{proof}  Let $\overleftarrow{Y}$ be a left enlargement of $Y$ and $f:\overleftarrow{Y} \to Y$ be an irreducible morphism.
 Consider $F(f):F(\overleftarrow{Y}) \to F(Y)$, where $F(\overleftarrow{Y})=X\oplus Y$ and $F(Y)=Y$. 
 
 Since $f$ is an irreducible morphism, by Proposition \ref{sec-ret-irr} we know that $F(f)$ is either a retraction or an irreducible morphism. 
First, suppose that $F(f)$ is irreducible. In this case, we show that $F(f)$ is not radical. In fact, otherwise if $F(f)$ is radical and we consider $i:Y \to X \oplus Y$ the inclusion morphism, since it is a section and $F(f)$ is radical-irreducible, then by Proposition \ref{BS3} we infer that the composite $F(f)i:Y \to Y$ is irreducible as well, which is a contradiction.
Thus, $F(f)$ is not radical and therefore is of the form $\left(\begin{array}{cc} c  & 1 \end{array} \right):X\oplus Y \to Y$ where $c:X\to Y$ is irreducible.

Now, assume that $F(f)$ is a retraction. Without loss of generality,  we write $F(f)=\left(\begin{array}{cc} c  & 1 \end{array} \right):X\oplus Y \to Y$ where $c:X\to Y$ is a morphism. In any case we have  a commutative diagram as follows:

 $${\tiny {\xymatrix { Z^{0} \ar[dd]^{0} \ar[rr]\ar[rr]&&  X^{1}\oplus Y^{1}  \ar[dd]^{\left(\begin{array}{cc} c_1  & 1 \end{array} \right)} \ar[rr]&& X^{2}\oplus Y^{2} \ar[dd]^{\left(\begin{array}{cc} c_2  & 1 \end{array} \right)}\ar[rr] &&
\dots \ar[rr]&&  X^{n} \oplus Y^{n} \ar[dd]^{\left(\begin{array}{cc} c_n  & 1 \end{array} \right)} &&\\ &&  &&    &&\\
0 \ar[rr] &&  Y^{1} \ar[rr] &&  Y^{2} \ar[rr] && \dots \ar[rr] &&  Y^{n}&&.}}}$$
\vspace{.1in}

Let $h:Y \to \overleftarrow{Y} $ be a morphism  given by $h^0=0$ and $h^i=\left(\begin{array}{c} 0 \\  1 \end{array} \right)$.  
It is not hard to check that $fh=Id_Y$ and that $f$ is a retraction, that contradicts our first assumption. Hence, there are not irreducible morphisms from $\overleftarrow{Y}$ to $Y$.
\end{proof}

Dually, one can prove that there are no irreducible morphisms in $\mathbf{C_{n+1}}(\mbox{proj} \,A)$ from an indecomposable complex $X$ in $\mathbf{C_{n}}(\mbox{proj} \,A)$  to a right enlargement of $X$.
\vspace{.1in}

Now, let $f:X \to Y$ be an irreducible morphism between indecomposable complexes in $\mathbf{C_{[0,n]}}(\mbox{proj} \,A)$. It is clear that if $X^0=Y^0=0$ then $X$ and $Y$ are in $\mathbf{C_{n}}(\mbox{proj} \,A)$ and $f$ is an irreducible morphism in $\mathbf{C_{n}}(\mbox{proj} \,A)$. Hence, we may consider the cases where $X^0\neq 0$ or $Y^0 \neq 0$. If $X^0=0$, then $X \in \mathbf{C_{n}}(\mbox{proj} \,A)$ and since $Y^0\neq 0$ then by Lemma {\ref{lema-previo}} we know that $Y$ is a left enlargement of $X$. Thus, $f:X \to \overleftarrow{X}$ is of the form (sec) where all entries are sections and $F(f)$ is a section of complexes in $\mathbf{C_{n}}(\mbox{proj} \,A)$.
 
 Now, if $Y^0=0$ then $X^0\neq 0$ and by Lemma \ref{noirred} the morphism $f:X \to Y$ is such that $F(f)$ is irreducible in $\mathbf{C_{n}}(\mbox{proj} \,A)$. In fact, if $F(f)$ is not an irreducible morphism then by Proposition \ref{sec-ret-irr} we infer that $F(f)$ is either a section or a retraction. This implies that $X=\overleftarrow{Y}$ and that $f:\overleftarrow{Y} \to Y$ is an irreducible morphism, which contradicts Lemma \ref{noirred}.  
 \vspace{.1in}

 Our next proposition deals with the last case where both $X^0$ and $Y^0$ are non-zero.
 
%\begin{prop}\label{clasif}
%Let $X$ and $Y$ be non stalk injective indecomposable complexes in $\mathbf{C_{[0,n]}} (\emph{proj} \,A)$ with $X^0\neq 0$ and $Y^0 \neq 0$. Let $f:X \to Y$ be an irreducible morphism in $\mathbf{C_{[0,n]}} (\emph{proj} \,A)$. Then  there exists an indecomposable complex $W$ in $\mathbf{C_{n}}(\emph{proj} \,A)$ such that $X$ and $Y$ are two different enlargements of $W$ and $f$ is a morphism between enlargements or there exists a radical-irreducible morphism $g$ in $\mathbf{C_{n}}(\emph{proj} \,A)$ such that $f^{0}$ is a retraction and $F(f)=g$.
%\end{prop} 

\begin{prop}\label{clasif}
Let $X$ and $Y$ be indecomposable complexes in $\mathbf{C_{[0,n]}} (\emph{proj} \,A)$ with $X^0\neq 0$ and $Y^0 \neq 0$ and $f:X \to Y$ be an irreducible morphism in $\mathbf{C_{[0,n]}} (\emph{proj} \,A)$. 

\begin{enumerate} 

\item[(a)] If $X$ and $Y$ are non stalk injective indecomposable complexes in $\mathbf{C_{[0,n]}} (\emph{proj} \,A)$ then  there exists an indecomposable complex $W$ in $\mathbf{C_{n}}(\emph{proj} \,A)$ such that $X$ and $Y$ are two different enlargements of $W$ and $f$ is a morphism between enlargements or there exists a radical-irreducible morphism $g$ in $\mathbf{C_{n}}(\emph{proj} \,A)$ such that $f^{0}$ is a retraction and $F(f)=g$.

\item[(b)] If $X$ or $Y$ are stalk injective complexes in $\mathbf{C_{[0,n]}} (\emph{proj} \,A)$ with $X^0\neq 0$, then $G(f)$ is  an irreducible morphism in $\mathbf{C_{[0,n-1]}} (\emph{proj} \,A)$. 
\end{enumerate}

\end{prop} 

\begin{proof}
(a) Let $f:X \to Y$ be an irreducible morphism. By Proposition \ref{sec-ret-irr} $ F(f):F(X)\to F(Y)$ is either a section, a retraction or an irreducible morphism in $\mathbf{C_{n}}(\mbox{proj} \,A)$. Assume that $F(f)$ is a section, then any direct summand  of $F(X)$ is also a direct summand of $F(Y)$. 

Since $X$ and $Y$ are not stalk injective complexes, we have that $F(X)\neq 0$ and  $F(Y)\neq 0$. Consider $W$ an indecomposable direct summand of $F(X)$, Therefore $W$ is a direct summand of $F(Y)$. Thus $X$ is a left enlargement of $W$, $X=(\overleftarrow{W})_1$ and $Y$ is another left enlargement of $W$, $Y=(\overleftarrow{W})_2$. Clearly in this case $f$ is a morphism between enlargements of $W$. Furthermore, we have a similar situation if $F(f)$ is a retraction. 
   
Now, assume that $F(f)$ is an irreducible morphism in $\mathbf{C_{n}}(\mbox{proj} \,A)$. If $F(f):F(X) \to F(Y)$ is not radical then by Proposition \ref{BS1} there exists a decomposition $F(f): W_1 \oplus X' \to W_1 \oplus Y'$ where $F(f)$ can be identified as the following matrix ${\left(\begin{array}{cc} 1 & 0 \\  0 & g  \end{array} \right)}$. Since $F(f)$ is not radical then $W_1 \neq 0$, and considering $W$ as any indecomposable direct summand of $W_1$ it follows that both $X$ and $Y$ are enlargements of $W$. Therefore, $f$ is a morphism between enlargements of the same complex.

Finally, let $F(f)$ be a radical-irreducible morphism in $\mathbf{C_{n}}(\mbox{proj} \,A)$. We shall show that if this is the case then $f^0$ is a retraction. In fact, consider the morphism $t:X^0 \to Y^1$ given by the composition $t=d ^0_yf^0=f^1d^0_X$.

The following commutative diagram:

$${\xymatrix { X^0 \ar[r]^{d ^0_X}\ar[d]^{1} & X^1 \ar[r]\ar[d]^{f^1} & \cdots \ar[r] & X^{n} \ar[d]^{f^n} \\  X^0 \ar[r]^{t}\ar[d]^{f^0} & Y^1 \ar[r]\ar[d]^{1} & \cdots \ar[r] & Y^{n} \ar[d]^{1} \\
Y^0 \ar[r]^{d ^0_y}& Y^1 \ar[r] & \cdots \ar[r] & Y^{n}}}$$ 
\vspace{.1in}

\noindent give us a factorization of $f$. Since $f$ is irreducible in $\mathbf{C_{[0,n]}}(\mbox{proj} \,A)$ either the first morphism is a section or the second is a retraction. In the former case, we get that $F(f)$ is a section, a contradiction to our previous assumption. In the latter case, $f^{0}$ is a retraction, proving the result.

(b) We only prove the statement when $Y$ is a stalk injective complex  since  the other case is trivial. Since $f:X \to Y$ is an irreducible morphism then by Proposition \ref{sec-ret-irr},   $G(f)$ is either a retraction, a section or irreducible. We shall show that $G(f)$ is not a retraction and hence is irreducible since clearly it is not a section. Assume  that is a retraction. Thus $Y$ is a direct summand of $G(X)$ but this  implies that $Y$ is a direct summand of $X$, a contradiction to the fact that $X$ is indecomposable. Hence $G(f)$ is an irreducible morphism.

\end{proof} 

A similar analysis can be done using right enlargements instead of left enlargements. Below, we state the dual of Proposition \ref{clasif}. 

\begin{prop}\label{dualclasif}
%Let $X$ and $Y$ be non stalk projective indecomposable complexes in $\mathbf{C_{[0,n]}}(\emph{proj} \,A)$ with $X^n\neq 0$ and $Y^n \neq 0$. Consider $f:X \to Y$ an irreducible morphism in $\mathbf{C_{[0,n]}}(\emph{proj} \,A)$. Then  there exists an indecomposable complex $W$ in $\mathbf{C_{[0,n-1]}}(\emph{proj} \,A)$ such that $X$ and $Y$ are two different right enlargements of $W$ and $f$ is a morphism between such enlargements or there exists a radical-irreducible morphism $g$ in $\mathbf{C_{[0,n-1]}}(\emph{proj} \,A)$ such that $f^{n}$ is a section and $G(f)=g.$ 

Let $X$ and $Y$ be indecomposable complexes in $\mathbf{C_{[0,n]}} (\emph{proj} \,A)$ with $X^n\neq 0$ and $Y^n \neq 0$ and $f:X \to Y$ be an irreducible morphism in $\mathbf{C_{[0,n]}} (\emph{proj} \,A)$. 

%\end{prop}
%\vspace{.1in}

\begin{enumerate} 

\item[(a)] If $X$ and $Y$ are non stalk projective indecomposable complexes in $\mathbf{C_{[0,n]}} (\emph{proj} \,A)$ then there exists an indecomposable complex $W$ in $\mathbf{C_{n}}(\emph{proj} \,A)$ such that $X$ and $Y$ are two different enlargements of $W$ and $f$ is a morphism between enlargements or there exists a radical-irreducible morphism $g$ in $\mathbf{C_{n}}(\emph{proj} \,A)$ such that $f^{n}$ is a section and $G(f)=g$.

\item[(b)] If $X$ or $Y$ are stalk projective complexes in $\mathbf{C_{[0,n]}} (\emph{proj} \,A)$ then $F(f)$ is an irreducible morphism in $\mathbf{C_{n}}(\emph{proj} \,A)$. 
\end{enumerate}

\end{prop}

\section{On indecomposability of diagonal complexes%in $\mathbf{C_{[0,n]}}(\mbox{proj} \,A)$ 
}

We study conditions under which there exists an enlargement of an indecomposable  complex $X$ in $\mathbf{C_{n}}(\mbox{proj} \,A)$.

\begin{prop} \label{existence} Let $X$ be an indecomposable  complex in $\mathbf{C_{n}}(\emph{proj} \,A)$. Assume that $X$ can be extended to the left and  $Z^{0} \in \emph{proj} \,A$ is such that ${\xymatrix { Z^{0}   \ar[r]^{\widetilde{d}_{X}^{0}}& X^{1} \ar[r]^{d_{X}^{1}} & \dots \ar[r]^{d_{X}^{n-1}} & X^{n} }}$ is a complex in $\mathbf{C_{[0,n]}}(\emph{proj} \,A)$ with $\widetilde{d}_{X}^{0}\neq0$.
Let $Y$ be a complex in $\mathbf{C_{n}}(\emph{proj} \,A)$ such that ${\xymatrix { Z^{0}   \ar[r]^{\widetilde{d}_{Y}^{0}}& Y^{1} \ar[r]^{d_{Y}^{1}} & \dots \ar[r]^{d_{Y}^{n-1}} & Y^{n} }}$ is an indecomposable  complex in $\mathbf{C_{[0,n]}}(\emph{proj} \,A)$. 

If $\emph{Hom}_{\tiny{K^{-,b}(\emph{proj}\,A)}}(X,Y)=0$ then there is a left enlargement $\overleftarrow{X}$ of $X$ in $\mathbf{C_{[0,n]}}(\emph{proj} \,A)$  of the form

%$${\xymatrix {&  X^{1} \ar[r]^{d_{X}^{1}} & \dots \ar[r]^{d_{X}^{n-1}} & X^{n}  \\ Z^{0}   \ar[ru]^{d_{X}^{0}} \ar[rd]^{d_{Y}^{0}} \\ & Y^{1} \ar[r]^{d_{Y}^{1}} & \dots \ar[r]^{d_{Y}^{n-1}} & Y^{n} }}$$

$${\xymatrix { Z^{0}  \ar[rr]^{\left(\begin{array}{cc} \widetilde{d}_{X}^{0} \\  \widetilde{d}_{Y}^{0}  \end{array} \right)}&&  X^{1}\oplus Y^{1}   \ar[rr]^{\left(\begin{array}{cc} d_{X}^{1} &0 \\  0 & d_{Y}^{1}  \end{array} \right)}&& X^{2}\oplus Y^{2}   \ar[rr]^{\left(\begin{array}{cc} d_{X}^{2} &0 \\  0 & d_{Y}^{2}  \end{array} \right)}&& \dots  \ar[rr]^{\left(\begin{array}{cc} d_{X}^{n-1} &0 \\  0 & d_{Y}^{n-1}  \end{array} \right)} &&  X^{n}\oplus Y^{n} }}.$$
\end{prop}

\begin{proof} Consider $h$ as the following  morphism:

 $${\xymatrix {\overleftarrow{Y}: \ar[d]^{h} & Z^{0} \ar[d]^{ \widetilde{d}_{X}^{0}} \ar[r]^{\widetilde{d}_{Y}^{0}}  & Y^{1}  \ar[d]^{0} \ar[r]^{d_{Y}^{1}}& Y^{2}  \ar[d]^{0} \ar[r]^{d_{Y}^{2}} & \dots \ar[r]& Y^{n} \ar[d]^{0 } \\
X: & X^{1} \ar[r]^{d_{X}^{1}} & X^{2} \ar[r] & \dots \ar[r] & X^{n}\ar[r] &0.}}$$
\vspace{.1in}

It is clear that $h$ is a non-zero morphism. Since $X$  is an  indecomposable complex,  then $ \widetilde{d}_{X}^{0}$ is not an isomorphism. Furthermore, { it follows from \cite[Corollary 1.4]{HZ}} that the cone of $h$, $C(h)$,  is indecomposable because $\mbox{Hom}(X,\overleftarrow{Y}[1])=0$ since $\mbox{Hom}(X,Y)=0$. Hence $C(h) = \overleftarrow{X}$ which is of the desired form.
\end{proof}
\vspace{.1in}

Next, we characterize the shape of $\overleftarrow{X}$  where in Proposition \ref{existence} we consider $X$ to be an indecomposable complex in $\mathbf{C_{n}}(\mbox{proj} \,A)$ and $Z^{0}$ an indecomposable projective $A$-module.
\vspace{.1in}

\begin{prop}\label{prop-indecomposable} Let $X_i$ for $i=1, \dots, m$ be indecomposable  complexes  in $\mathbf{C_{n}}(\emph{proj} \,A)$ and $Z^{0}$ be an indecomposable projective  $A$-module. 
Then the diagonal complex
$${\xymatrix {Z:  Z^{0}  \ar[r] & X_1 \oplus \dots \oplus X_{m}  }}$$ is not indecomposable if and only if there exists  $i=1, \dots, m$ such that $X_i$ in $\mathbf{C_{[0,n]}}(\emph{proj} \,A)$ is a direct summand of $Z$.
\end{prop}
%\begin{proof} Let $Z$  be {\color{purple}an indecomposable?}  direct summand of $W$ in $\mathbf{C_{[0,n]}}(\mbox{proj} \,A)$. Then $F(Z)$ is a direct summand of $F(W)=\oplus_{i=1}^{m} W_{i}$ in $\mathbf{C_{n}}(\mbox{proj} \,A)$. Therefore $F(Z)=\oplus_{j=1}^{r} W_{i_{j}}$.

%In case $F(Z)=0$, that is, $r=0$, then $Z=0$ or $Z=W^{0}\rightarrow 0 \rightarrow0 \dots \rightarrow 0$. 

%The first case is clearly not possible.  Now, if the latter case holds, then $W=Z\oplus Z^{\prime}$ where $Z^{\prime}= \oplus_{i=1}^{m} W_{i}$.  Hence, for all $j=1, \dots, m$ we see that $W_{j}$ is a direct summand of $W$.

%In case $F(Z)=\oplus_{i=1}^{m} W_{i}$, that is $r=m$, then $Z= \oplus_{i=1}^{m} W_{i}$ or $Z=W$. 
%The first case is not possible because $Z$ is indecomposable.  In the latter case, we have that $Z$ is a trivial summand. 

%Without loss of generality, we may assume that $F(Z)=\oplus_{j=1}^{r} W_{i_{j}}$ with $0 <r< m$.   Since $Z$ is indecomposable then $r=1$ or $Z^{0}=W^{0}$.

%In case $r=1$ then $Z=0W_{j}$ for some $j$ or $W^{0} \simeq Z^{0}$.

%If $W^{0} \simeq Z^{0}$, then as we analyzed before, there is a complex $Z^{\prime}$ such that $W=Z\oplus Z^{\prime}$ and $Z^{\prime}= \oplus_{i\neq i_{j}}^{m} W_{i}$.  Hence, there is a positive integer $k$ where $W_{k}$ is a direct summand of $W$.
%\end{proof}

%{\color{purple} OTRA FORMA DE ESCRIBIR LA PRUEBA MAS CORTA }

\begin{proof}
 Let $W$  be an indecomposable direct summand of $Z$ in $\mathbf{C_{[0,n]}}(\mbox{proj} \,A)$. Then $F(W)$ is a direct summand of $F(Z)=\oplus_{i=1}^{m} X_{i}$ in $\mathbf{C_{n}}(\mbox{proj} \,A)$. Therefore $F(W)=\oplus_{j=1}^{r} X_{i_{j}}$.

First, assume that $W^0=0$, then $r=1$ since $W$ is indecomposable. Therefore, $W=X_i$ for some $i$ and the claim follows. 

 Assume now that $W^0=Z^0$. Then there is a complex $W^{\prime}$ such that $Z=W\oplus W^{\prime}$ where $W^{\prime}= \oplus_{i\neq i_{j}}  X_{i}$. If $r=m$ then  $Z$ is indecomposable because in this case $Z=W$. Otherwise,  there is a positive integer $k$ where $X_{k}$ is a direct summand of $Z$.
\end{proof}

Finally, we analyze the conditions under which the complex ${\xymatrix { 0   \ar[r]^{d_{X}^{0}}& X^{1} \ar[r]^{d_{X}^{1}} & \dots \ar[r]^{d_{X}^{n-1}} & X^{n}}}$ is direct summand of a complex $Z$ in $\mathbf{C_{[0,n]}}(\mbox{proj} \,A)$.

\begin{lem}\label{prop-summads} Let $X$ be indecomposable  complex in $\mathbf{C_{n}}(\emph{proj} \,A)$ and  
$${\xymatrix {Z:  Z^{0}  \ar[rr]^{\left(\begin{array}{cc} \widetilde{d}_{X}^{0} \\  \widetilde{d}_{Y}^{0}  \end{array} \right)}&&  X^{1}\oplus Y^{1}   \ar[rr]^{\left(\begin{array}{cc} d_{X}^{1} &0 \\  0 & d_{Y}^{1}  \end{array} \right)}&& X^{2}\oplus Y^{2}   \ar[rr]^{\left(\begin{array}{cc} d_{X}^{2} &0 \\  0 & d_{Y}^{2}  \end{array} \right)}&& \dots  \ar[rr]^{\left(\begin{array}{cc} d_{X}^{n-1} &0 \\  0 & d_{Y}^{n-1}  \end{array} \right)} &&  X^{n}\oplus Y^{n} }}$$ be a complex in $\mathbf{C_{[0,n]}}(\emph{proj} \,A)$. 
\begin{enumerate}
    \item[(i)] The complex $X$ is a direct summand of $Z$ if and only if there exists $g:Y \rightarrow X$ in $\mathbf{C_{n}}(\emph{proj} \,A)$ such that $-g^{1}\widetilde{d}_{Y}^{0}=\widetilde{d}_{X}^{0}$.
    \item[(ii)] If $X$ is a direct summand of $Z$, if and only if $Z$ is isomorphic to another complex of the same shape but with the following differentials: $${\xymatrix {\widetilde{Z}:  Z^{0}  \ar[rr]^{\left(\begin{array}{cc} {0} \\  \widetilde{d}_{Y}^{0}  \end{array} \right)}&&  X^{1}\oplus Y^{1}   \ar[rr]^{\left(\begin{array}{cc} d_{X}^{1} &0 \\  0 & d_{Y}^{1}  \end{array} \right)}&& X^{2}\oplus Y^{2}   \ar[rr]^{\left(\begin{array}{cc} d_{X}^{2} &0 \\  0 & d_{Y}^{2}  \end{array} \right)}&& \dots  \ar[rr]^{\left(\begin{array}{cc} d_{X}^{n-1} &0 \\  0 & d_{Y}^{n-1}  \end{array} \right)} &&  X^{n}\oplus Y^{n}}}.$$
\end{enumerate}
\end{lem}

\begin{proof} (i) Assume that ${\xymatrix { 0   \ar[r]^{d_{X}^{0}}& X^{1} \ar[r]^{d_{X}^{1}} & \dots \ar[r]^{d_{X}^{n-1}} & X^{n}}}$ is a direct summand of $Z$. Then there is a section $h:  ({\xymatrix { 0   \ar[r]^{d_{X}^{0}}& X^{1} \ar[r]^{d_{X}^{1}} & \dots \ar[r]^{d_{X}^{n-1}} & X^{n}}}) \longrightarrow Z$. Hence, there exists a morphism $t: Z\longrightarrow ( {\xymatrix { 0   \ar[r]^{d_{X}^{0}}& X^{1} \ar[r]^{d_{X}^{1}} & \dots \ar[r]^{d_{X}^{n-1}} & X^{n}}})$ such that $th= 1_{X}$. Moreover $t$ is such that $t^{0}=0$ and $t^{i}={\left(\begin{array}{cc} 1 & g^{i}  \end{array} \right)}$ where $g^{i}:Y^{i}\rightarrow X^{i}$. Since $d_{X}^{i} {\left(\begin{array}{cc} 1 & g^{i}  \end{array} \right)}= {\left(\begin{array}{cc} 1 & g^{i+1}  \end{array} \right)}{\left(\begin{array}{cc} d_{X}^{i} &0 \\  0 & d_{Y}^{i}  \end{array} \right)}$ then  ${\left(\begin{array}{cc} d_{X}^{i} &  d_{X}^{i} g^{i}  \end{array} \right)}= {\left(\begin{array}{cc} d_{X}^{i} & g^{i+1}d_{Y}^{i}  \end{array} \right)}$ proving that $g:Y \rightarrow X$ is a morphism of complexes.

On the other hand, the morphism $g:Y \rightarrow X$  satisfies ${\left(\begin{array}{cc} 1 & g^{1}  \end{array} \right)} {\left(\begin{array}{c} \widetilde{d}_{X}^{0}\\ \widetilde{d}_{Y}^{0}  \end{array} \right)}=0$. Therefore, we have $-g^{1}\widetilde{d}_{Y}^{0}=\widetilde{d}_{X}^{0}$.

The converse is clear, since under our assumptions the morphism $t=\{(1,g_i)\}_{i=1}^n: Z\rightarrow X$ is a retraction. 

(ii) It follows from Statement (i) that there exists a morphism $g:Y \rightarrow X$ in $\mathbf{C_{n}}(\mbox{proj} \,A)$ such that $-g^{1}\widetilde{d}_{Y}^{0}=\widetilde{d}_{X}^{0}$. Then we can construct an isomorphism of complexes as follows: 
$${\xymatrix {Z^{0} \ar[rr]^{\left(\begin{array}{cc} \widetilde{d}_{X}^{0} \\  \widetilde{d}_{Y}^{0}  \end{array} \right)} \ar[d]^{1} && X^{1} \oplus Y^{1} \ar[rr] \ar[d]^{\left(\begin{array}{cc} {1} & g^{1} \\  0  & {1}  \end{array} \right)} && X^{2}\oplus Y^{2}\ar[rr] \ar[d]^{\left(\begin{array}{cc} {1} & g^{2} \\  0  & {1}  \end{array} \right)}&& \dots \ar[rr] && \ar[d] X^{n}\oplus Y^{n} \\ Z^{0}   \ar[rr]^{\left(\begin{array}{cc} {0} \\  \widetilde{d}_{Y}^{0}  \end{array} \right)}&&  X^{1}\oplus Y^{1} \ar[rr]&& X^{2}\oplus Y^{2}   \ar[rr]&& \dots  \ar[rr] &&  X^{n}\oplus Y^{n}  }}$$
proving the result.
\end{proof}

 As an immediate consequence of Proposition \ref{prop-indecomposable} and Lemma \ref{prop-summads} we get the following corollary.

%\begin{coro}\label{coro-indecomposable} Let $X, Y$ be indecomposable  complexes  in $\mathbf{C_{n}}(\emph{proj} \,A)$ and $Z^{0}$ be an indecomposable projective  $A$-module. 
%Then the complex
%$${\xymatrix {W:  Z^{0}  \ar[rr]^{\left(\begin{array}{cc} \widetilde{d}_{X}^{0} \\  \widetilde{d}_{Y}^{0}  \end{array} \right)}&&  X^{1}\oplus Y^{1}   \ar[rr]^{\left(\begin{array}{cc} d_{X}^{1} &0 \\  0 & d_{Y}^{1}  \end{array} \right)}&& X^{2}\oplus Y^{2}   \ar[rr]^{\left(\begin{array}{cc} d_{X}^{2} &0 \\  0 & d_{Y}^{2}  \end{array} \right)}&& \dots  \ar[rr]^{\left(\begin{array}{cc} d_{X}^{n-1} &0 \\  0 & d_{Y}^{n-1}  \end{array} \right)} &&  X^{n}\oplus Y^{n} }}$$ is decomposable  if and only if there are not a morphisms $g: Y \rightarrow X$ and  $h: X \rightarrow Y$ in  $\mathbf{C_{n}}(\emph{proj} \,A)$ such that  $-g^{1}\widetilde{d}_{Y}^{0}=\widetilde{d}_{X}^{0}$ and  $-h^{1}\widetilde{d}_{X}^{0}=\widetilde{d}_{Y}^{0}$.
%\end{coro} 
%{\color{purple}INTENTO DE REESCRIBIR EL COROLARIO 
\begin{coro}\label{coro-indecomposable}  Let 
${\xymatrix {Z:  Z^{0}  \ar[r] & X_1 \oplus \dots \oplus X_{m}  }}$ be a complex, where $X_i\in \mathbf{C_{n}}(\emph{proj} \,A)$ is indecomposable for each $i$ and and $Z^{0}$ is an indecomposable projective  $A$-module. 
 $Z$ is decomposable  if and only if for some $i$ there exists a morphism $g_i:\hat{X_i}= \oplus_{j\neq i} X_j\rightarrow X_i$  such that  $-g_i^{1}\widetilde{d}_{\hat{X_i}}^{0}=\widetilde{d}_{X_i}^{0}$.
\end{coro}

\section{Description of an enlargement}

In this section we study which projective $A$-modules are involved in  $Z^{0}$ to construct an enlargement $\overleftarrow{X}$ of a complex $X$ in $\mathbf{C_{n}}(\mbox{proj} \,A)$. %via irreducible morphisms. %Furthermore, we show how is the shape of an enlargement.

\begin{prop}\label{prop-summads2} Let $X$ be an indecomposable  complex in $\mathbf{C_{n}}(\emph{proj} \,A)$ and $f:X \rightarrow \overleftarrow{X}$
be an irreducible morphism in $\mathbf{C_{[0,n]}}(\emph{proj} \,A)$. Assume that $\overleftarrow{X}$ is as follows $${\xymatrix {\overleftarrow{X}:
Z^{0}  \ar[rr]^{\left(\begin{array}{cc} \widetilde{d}_{X}^{0} \\  \widetilde{d}_{Y}^{0}  \end{array} \right)}&&  X^{1}\oplus Y^{1}
 \ar[rr]^{\left(\begin{array}{cc} d_{X}^{1} &0 \\  0 & d_{Y}^{1}  \end{array} \right)}&& X^{2}\oplus Y^{2}   \ar[rr]^{\left(\begin{array}{cc} d_{X}^{2} &0 \\
   0 & d_{Y}^{2}  \end{array} \right)}&& \dots  \ar[rr]^{\left(\begin{array}{cc} d_{X}^{n-1} &0 \\  0 & d_{Y}^{n-1}  \end{array} \right)} &&  X^{n}\oplus Y^{n}}}.$$ Then $Z^{0}$ is a direct summand of the projective cover of $\emph{soc}(\emph{Ker}\, {d}_{X}^{1})$.
\end{prop}
\begin{proof} %Consider
%$${\xymatrix {\overleftarrow{X}:  Z^{0}  \ar[rr]^{\left(\begin{array}{cc} \widetilde{d}_{X}^{0} \\  \widetilde{d}_{Z}^{0}  \end{array} \right)}&&  X^{1}\oplus Z^{1}   \ar[rr]^{\left(\begin{array}{cc} d_{X}^{1} &0 \\  0 & d_{Z}^{1}  \end{array} \right)}&& X^{2}\oplus Z^{2}   \ar[rr]^{\left(\begin{array}{cc} d_{X}^{2} &0 \\  0 & d_{Z}^{2}  \end{array} \right)}&& \dots  \ar[rr]^{\left(\begin{array}{cc} d_{X}^{n-1} &0 \\  0 & d_{Z}^{n-1}  \end{array} \right)} &&  X^{n}\oplus Z^{n} }}.$$

Since $d_{X}^{1} \widetilde{d}_{X}^{0}=0$ then $\widetilde{d}_{X}^{0}$ factors through $\mbox{Ker}\, {d}_{X}^{1}$. Moreover, $Im \, \widetilde{d}_{X}^{0}\subset \mbox{Ker}\, {d}_{X}^{1}.$ Then  $\widetilde{d}_{X}^{0}= \iota h^{1}$ where $h^{1}:Z^{0} \rightarrow \mbox{Ker}\, {d}_{X}^{1}$ and $\iota: \mbox{Ker}\, {d}_{X}^{1} \rightarrow X^{1}$. % is a monomorphism since $\mbox{Ker}\, {d}_{X}^{1}= Im \, \widetilde{d}_{X}^{0}.$

Now, we can consider the non-zero morphism $Z^{0} \twoheadrightarrow Im \,h^{1} \hookrightarrow \mbox{Ker}\, {d}_{X}^{1}$ where $$0\neq \mbox{soc}(Im\, h^{1})  \subset \mbox{soc}(\mbox{Ker}\, {d}_{X}^{1}).$$

Let $\widetilde{P}$ be the projective cover of $\mbox{soc}(\mbox{Ker}\, {d}_{X}^{1}).$ We can consider the non-zero  composition of morphisms  $\widehat{\iota} \pi \widetilde{\pi}$ where $\widetilde{\pi}: \widetilde{P} \rightarrow \mbox{soc}(\mbox{Ker}\, {d}_{X}^{1})$, $\pi: \mbox{soc}(\mbox{Ker}\, {d}_{X}^{1})\rightarrow \mbox{soc}(Im\, h^{1})$ and $\widehat{\iota}: \mbox{soc}(Im\, h^{1}) \rightarrow Im\, h^{1}$.

On the other hand, applying $\mbox{Hom}(\widetilde{P}, -)$ to the epimorphism  $Z^{0} \rightarrow Im\, h^{1}$ we get the following epimorphism: $$\mbox{Hom}(\widetilde{P}, Z^{0}) \rightarrow \mbox{Hom}(\widetilde{P}, Im\, h^{1}).$$

\noindent Therefore, there is a morphism $t: \widetilde{P} \rightarrow  Z^{0}$ such that  $h^{1}t= \widehat{\iota} \pi \widetilde{\pi}$.

Now, we can construct the following  complex
$${\xymatrix { W: \widetilde{P}  \ar[rr]^{\left(\begin{array}{cc} \iota \widehat{\iota} \pi \widetilde{\pi} \\  \widetilde{d}_{Y}^{0} t \end{array} \right)}&&  X^{1}\oplus Y^{1}   \ar[rr]^{\left(\begin{array}{cc} d_{X}^{1} &0 \\  0 & d_{Y}^{1}  \end{array} \right)}&& X^{2}\oplus Y^{2}   \ar[rr]^{\left(\begin{array}{cc} d_{X}^{2} &0 \\  0 & d_{Y}^{2}  \end{array} \right)}&& \dots  \ar[rr]^{\left(\begin{array}{cc} d_{X}^{n-1} &0 \\  0 & d_{Y}^{n-1}  \end{array} \right)} &&  X^{n}\oplus Y^{n} }}.$$
It is not hard to see that the irreducible morphism $f:X \rightarrow \overleftarrow{X}$ factors through such a complex. %We illustrate the situation with the following diagram:
Then $f=g_2 g_1$ where $g_1: X\rightarrow W$ and $g_2: W \rightarrow \overleftarrow{X}$.

Now, if $g_1$ is a section then $X$ is a direct summand of $W$. Moreover,  by Lemma \ref{prop-summads}, we know that $W$ is isomorphic to a complex $W^{\prime}$ of the form:  
$${\xymatrix { W^{\prime}: \widetilde{P}  \ar[rr]^{\left(\begin{array}{cc} 0 \\  \widetilde{d}_{Y}^{0} t \end{array} \right)}&&  X^{1}\oplus Y^{1}   \ar[rr]^{\left(\begin{array}{cc} d_{X}^{1} &0 \\  0 & d_{Y}^{1}  \end{array} \right)}&& X^{2}\oplus Y^{2}   \ar[rr]^{\left(\begin{array}{cc} d_{X}^{2} &0 \\  0 & d_{Y}^{2}  \end{array} \right)}&& \dots  \ar[rr]^{\left(\begin{array}{cc} d_{X}^{n-1} &0 \\  0 & d_{Y}^{n-1}  \end{array} \right)} &&  X^{n}\oplus Y^{n} }}.$$
Therefore, we have an isomorphism of complexes  $\delta:W\rightarrow W^{\prime}$. Hence we can write $f=g_2g_1=(g_2 \delta^{-1})(\delta g_1)$. We denote by $g_2'=g_2 \delta^{-1}$ and by $g_1'= \delta g_1$. Furthermore, since we are assuming that $g_1$ is a section then either so is $g_1'$.

Thus, $d_{X}^{0} t \delta_0^{-1}=0$ a contradiction since $d_{X}^{0} t = i \widehat{i} \pi \widetilde{\pi}\neq 0$ and $ \delta^{-1}$ is an isomorphism. Hence, we prove that  $g_1$ is not a section.

Since $f$ is irreducible, then we have that $g_2$ is a retraction. Therefore $Z^{0}$ is a direct summand of $\widetilde{P}$, as we wish to prove.
\end{proof}

Now, we are in position to prove the main result of this section.

\begin{teo}
Let $X$ be an indecomposable  complex in $\mathbf{C_{n}}(\emph{proj} \,A)$ and $f:X \rightarrow Z$
be a left almost split morphism in $\mathbf{C_{[0,n]}}(\emph{proj} \,A)$. Then any indecomposable direct summand of the projective cover of $\emph{soc}(\emph{Ker}\, {d}_{X}^{1})$ is a direct summand of $Z^0$. Moreover, if $X$ is not injective then $(\tau^{-1}X)^0$ is in $\mbox{add}(P_S)$ %a direct sum of all the projective involve in
where $P_S$ is the projective cover of $\emph{soc}(\emph{Ker}\, {d}_{X}^{1})$.
% with some of them possibly repeat?. 
\end{teo}
\begin{proof}
Let $P$ be an indecomposable direct summand of the projective cover of $\mbox{soc}(\mbox{Ker}\, {d}_{X}^{1})$. Then $W:{\xymatrix { P   \ar[r]^{t}& X^{1} \ar[r]^{d_{X}^{1}} & \dots \ar[r]^{d_{X}^{n-1}} & X^{n}}}$ is an indecomposable complex, where the differential $t$ is the usual morphism that sends the top of $P$ to the socle  of $P$.  Note that the image of $t$ is simple. Moreover, if $ st\neq0$ for any non-zero morphism $s$ then $s$ is an epimorphism. 

Consider the morphism $\hat{f}:X\rightarrow W$ defined as follows: 
$${\xymatrix {0\ar[rr] \ar[d]^{0} && X^{1} \ar[rr] \ar[d]^1 && X^{2} \ar[rr]\ar[d]^{1}&& \dots \ar[rr] && \ar[d]^1 X^{n} \\ P   \ar[rr]^t&&  X^{1}\ar[rr]&& X^{2}  \ar[rr]&& \dots  \ar[rr] &&  X^{n}  }}.$$

Observe that the morphism $\hat{f}$ is not a section. Therefore, it factors through $f$, since $f$ is a left almost split morphism: 
$${\xymatrix { 0 \ar[d]
\ar[rr] && X^{1} \ar[d]^{g_1}  \ar[rr] && X^{2}\ar[rr]\ar[d]^{g_2} && \dots \ar[rr] && \ar[d]^{g_n} X^{n} \\ Z^{0} \ar[d]^{h_0}  \ar[rr]^{d^0_Z}&&  Z^{1} \ar[d]^{h_1}  \ar[rr]^{d^1_Z}&&  Z^{2}  \ar[d]^{h_2} \ar[rr]&& \dots  \ar[rr] &&   Z^{n} \ar[d]^{h_n} \\
P \ar[rr]^t && X^{1}  \ar[rr] && X^{2} \ar[rr] && \dots \ar[rr] && X^{n}}}.$$

Then $h_ig_i=1$ for all $i=1, \dots, n$ and $h_0t=h_1d^0_Z$. If $h_1d^0_Z=0$ then $X$ is a direct summand of $Z$, since the following morphism of complexes is a section: 
$${\xymatrix { 0 \ar[d]
\ar[rr] && X^{1} \ar[d]^{g_1}  \ar[rr] && X^{2}\ar[rr]\ar[d]^{g_2} && \dots \ar[rr] && \ar[d]^{g_n} X^{n} \\ Z^{0}   \ar[rr]^{d^0_Z}&&  Z^{1}  \ar[rr]^{d^1_Z}&&  Z^{2}  \ar[rr]&& \dots  \ar[rr] &&   Z^{n}  \\
}}.$$
Therefore, we get a contradiction to the fact that $f$ is irreducible and there are no irreducible morphisms from $X$ to $X$ in $\mathbf{C_{[0,n]}}(\mbox{proj} \,A)$. Then, $h_1d^0_Z\neq 0$ and thus $h_0t \neq 0$. Hence $h_0$ is an epimorphism. Therefore, $P$ is a direct summand of $Z^0$ as we wish to prove. 
\end{proof}

\begin{coro}
Let $X$ be indecomposable  complex in $\mathbf{C_{n}}(\emph{proj} \,A)$ and $f:X \rightarrow Z$
be a left almost split morphism in  $\mathbf{C_{n}}(\emph{proj} \,A)$. If $f:X \rightarrow Z$
is a left almost split morphism in  $\mathbf{C_{[0,n]}}(\emph{proj} \,A)$ then $\emph{Ker}d_x^1=0.$
\end{coro}

%In case that a morphism $f:X \rightarrow Y$  in $\mathbf{C_{n}}(\mbox{proj} \,A)$ is not irreducible and moreover
%$f^{1} \mid _{\mbox{soc}(\mbox{Ker}\, {d}_{X}^{1})}$ is not a monomorphism then by \cite[Theorem 3.9]{CGP6} we know that $f$ factors through an
%extension of $X$. Therefore if we assume that for a non-irreducible morphism we have that $f^{1} \mid _{\mbox{soc}(\mbox{Ker}\, {d}_{X}^{1})}$ is a
%monomorphism then it means that $f$ factors through a proper enlargement.
\vspace{.1in}

Next, we study the possible differentials $d_{Z}^{0}$ for a particular left enlargement $Z$, where $Z^{0}$ is a direct sum of two indecomposable projective $A$-modules. 

We observe that these complexes appear in a left almost split morphism of $\mathbf{C_{[0,n]}}(\mbox{proj} \,A)$.

\begin{prop}\label{useful2} Let  $X$ be an indecomposable complex in $\mathbf{C_{n}}(\emph{proj} \,A)$. Assume that there are not morphisms from an extension of $X$ to an indecomposable  complex $Z$, where $Z$ is as follows:

$${\xymatrix {  Z^{0}_{0} \oplus Z_{0}^{1} \ar[rr]^{\left(\begin{array}{cc}a & c \\  b & d  \end{array} \right)}&&  X^{1}\oplus W^{1}   \ar[rr]^{\left(\begin{array}{cc} d_{X}^{1} &0 \\  0 & d_{W}^{1}  \end{array} \right)}&& X^{2}\oplus W^{2}   \ar[rr]^{\left(\begin{array}{cc} d_{X}^{2} &0 \\  0 & d_{W}^{2}  \end{array} \right)}&& \dots  \ar[rr]^{\left(\begin{array}{cc} d_{X}^{n-1} &0 \\  0 & d_{W}^{n-1}  \end{array} \right)} &&  X^{n}\oplus W^{n}}} $$ with $Z^{0}= Z_{0}^{0} \oplus Z_{0}^{1}$ and where $ Z_{0}^{0}$ and $ Z_{0}^{1}$ are indecomposable projective $A$-modules.

Then $b\neq 0$ and $d\neq 0$. That is, the possibilities for $d_Z^0$ are the matrices ${\left(\begin{array}{cc}a & 0 \\  b & d  \end{array} \right)}$, 
${\left(\begin{array}{cc}a & c \\  b & d  \end{array} \right)}$ or ${\left(\begin{array}{cc}0 & c \\  b & d  \end{array} \right)}$. In all cases stated above for the possible differentials $d_{Z}^{0}$, the complexes $Z$ are not  isomorphic to each other.
%\begin{enumerate}

%\item If $a\neq 0$ then $b\neq 0$. Moreover,
%\begin{enumerate}
%\item if $c=0$ then $d\neq 0$.
%\item if $c\neq 0$ then $d\neq 0$.
%\end{enumerate}

%\item If $a=0$ then $b\neq 0$, $c\neq 0$ and  $d\neq 0$.
%\end{enumerate}

\end{prop}

\begin{proof} If $a\neq 0$ then $b\neq 0$. In fact, the complex  $Z_{0}^{0}\rightarrow  X^1 \rightarrow \dots  \rightarrow X^n$ is an extension of $X$, since $Z_{0}^{0}$ is an indecomposable projective $A$-module. Therefore, if we assume that $b=0$ then  there is a non-zero morphism from an extension of $X$ to $Z$ as we illustrate in  the following diagram:

$${\xymatrix {   Z_{0}^{0} \oplus  Z_{0}^{1} \ar[rr]^{\left(\begin{array}{cc}a & c \\  0 & d  \end{array} \right)}&&  X^{1}\oplus W^{1}   \ar[rr]^{\left(\begin{array}{cc} d_{X}^{1} &0 \\  0 & d_{W}^{1}  \end{array} \right)}&& X^{2}\oplus W^{2}   \ar[rr]^{\left(\begin{array}{cc} d_{X}^{2} &0 \\  0 & d_{W}^{2}  \end{array} \right)}&& \dots \\ && && \\   Z_{0}^{0}  \ar[rr]^{a} \ar[uu]^{\left(\begin{array}{c}1 \\  0  \end{array} \right)}  &&  X^{1} \ar[rr]^{d_{X}^{1}} \ar[uu]^{\left(\begin{array}{c}1 \\  0  \end{array} \right)} && X^{2}  \ar[rr]^{d_{X}^{2}}  \ar[uu]^{\left(\begin{array}{c}1 \\  0  \end{array} \right)} && \dots }}$$

\noindent which clearly is a morphism of complexes.

Moreover, $c\neq 0$ or $d\neq 0$. In fact otherwise if both are zero-morphisms then the complex $Z$ is not indecomposable. Therefore, if $c=0$ then $d\neq 0$.

Now, assume that $c\neq 0$. Then $d\neq 0$. Indeed, the complex $Z_{0}^{1} \rightarrow  X^1 \rightarrow \dots  \rightarrow X^n$ is an extension of $X$, since $Z_{0}^{1}$ is an indecomposable projective $A$-module. Therefore,  there is a non-zero morphism of complexes from an extension of $X$ to $Z$ as we illustrate below:

$${\xymatrix {  Z^{0}_{0} \oplus  Z_{0}^{1} \ar[rr]^{\left(\begin{array}{cc}a & c \\  b & d  \end{array} \right)}&&  X^{1}\oplus W^{1}   \ar[rr]^{\left(\begin{array}{cc} d_{X}^{1} &0 \\  0 & d_{W}^{1}  \end{array} \right)}&& X^{2}\oplus W^{2}   \ar[rr]^{\left(\begin{array}{cc} d_{X}^{2} &0 \\  0 & d_{W}^{2}  \end{array} \right)}&& \dots \\ && && \\   Z_{0}^{1}  \ar[rr]^{c} \ar[uu]^{\left(\begin{array}{c}0 \\  1  \end{array} \right)}  &&  X^{1} \ar[rr]^{d_{X}^{1}} \ar[uu]^{\left(\begin{array}{c}1 \\  0  \end{array} \right)} && X^{2}  \ar[rr]^{d_{X}^{2}}  \ar[uu]^{\left(\begin{array}{c}1 \\  0  \end{array} \right)} && \dots }}$$

\noindent  where the first square commutes if and only if $d=0$.

In conclusion, $d_{Z}^{0}={\left(\begin{array}{cc}a & 0 \\  b & d  \end{array} \right)}$ or
$d_{Z}^{0}={\left(\begin{array}{cc}a & c \\  b & d  \end{array} \right)}$.

Now, if $a=0$ then $b\neq 0$, since otherwise the complex $Z$ is not indecomposable. Moreover, for the same reason we can see that $c\neq 0$ and hence also $d\neq 0$. In conclusion, in this case we get that $d_{Z}^{0}={\left(\begin{array}{cc}0 & c \\  b & d  \end{array} \right)}$.

We leave to the reader to verify that all possible complexes $Z$ are not isomorphic between each other.
\end{proof}

\section{ Irreducible morphism between enlargements}

In this section we shall analyze how are the entries of an irreducible morphism $f$ in $\mathbf{C_{[0, n]}}(\mbox{proj} \,A)$ between enlargements from a complex
$X \in \mathbf{C_{n}}(\mbox{proj} \,A)$. We shall give a method to decide the form of such irreducible morphism.
For this purpose it will be fundamental Proposition \ref{GM}.
\vspace{.1in}

We start with the following proposition.

\begin{prop}\label{irred} Let $A$ be an artin algebra. Let $X$ be an indecomposable complex in $\mathbf{C_{n}}(\emph{proj} \,A)$,
$(\overleftarrow{X})_1 : Z_0 \rightarrow  X \oplus W $ and
$(\overleftarrow{X})_2 : Z'_{0} \rightarrow X \oplus W'$ be enlargements of $X$, where $W$ and $W'$ are  non-zero complexes in $\mathbf{C_{n}}(\emph{proj} \,A)$. Moreover, $W$ and $W'$ have  non-zero non-common direct summands.  

Let $f
%=\left(
%     \begin{array}{cc}
%        f^{0}, & F(f) \\
%     \end{array}
%   \right)
: (\overleftarrow{X})_1 \rightarrow (\overleftarrow{X})_2$ be an irreducible morphism in $\mathbf{C_{[0, n]}}(\emph{proj} \,A)$. Then  $F(f) \simeq \left(
      \begin{array}{cc}
        1 & 0 \\
        0 & g \\
      \end{array}
    \right)$ with $g: W \rightarrow W'$ a radical-irreducible morphism.
\end{prop}

\begin{proof} Assume that $f$ is irreducible. By Proposition \ref{sec-ret-irr} we have that $F(f)$ is a section, a retraction or an irreducible morphism.
If $F(f)$ is a section then $W$ is a direct summand of $W'$ and if $F(f)$ is a retraction then $W'$ is a direct summand of $W$, a contradiction. Then $F(f)$ is  irreducible.

We claim that the morphism $F(f):X \oplus W \rightarrow  X \oplus W'$ is not a radical morphism. In fact, otherwise, if $F(f)$  is a radical-irreducible morphism then by Proposition \ref{BS3} the composition of $F(f)$ with any section is also irreducible. Consider the section given by the inclusion $i:X \to X\oplus W$, we have that $F(f)i:X\to X \oplus W'$ is also irreducible. We observe that if  $F(f)i$ is not a radical morphism, then since $X$ is indecomposable there exists $\Bar{W}$ an indecomposable direct summand of $X\oplus W'$ and a retraction $\pi:W \to \Bar{W}$ such that the composition $\pi F(f)i$ is an isomorphism, a contradiction to $F(f)$ being a radical morphism. Hence $F(f)i$ is a radical-irreducible morphism, therefore again by  Proposition \ref{BS3} the composition $X \stackrel{F(f)i}\to X\oplus W' \to X$ is irreducible. A contradiction since by  \cite[Theorem 2.19]{BS} there are not irreducible morphisms from $X$ to $X$ when $X$ is an indecomposable complex. 

Hence $F(f)$ is not a radical morphism. Therefore by Proposition \ref{BS1} $F(f)$ is isomorphic to a morphism of the form: 
$$F(f)=\left(
      \begin{array}{cc}
        1 & 0 \\
        0 & g \\
      \end{array}
    \right) : X \oplus W \rightarrow  X \oplus W'$$

\noindent with $g:W \rightarrow W'$ a radical morphism. Since $f$ is irreducible then by  Proposition \ref{BS2}, we have that $g$ is irreducible.
\end{proof}

The converse does not hold as we show in the next example. 

\begin{ej}
    \emph{ Let $A$ be the algebra given by the quiver:}
{\begin{displaymath}
\xymatrix @R=0.1cm  @C=1cm{&&  3 \ar[r]^{\alpha_3}&4 \\  1 \ar[r]^{\alpha_1} & 2  \ar[rd]_{\beta_1}\ar[ru]^{\alpha_2}&  \\
 & &   5\ar[r]_{\beta_2}&6}
\end{displaymath}}
 \emph{with the relations $\alpha_3\alpha_2\alpha_1=\beta_2 \beta_1 \alpha_1 =0$.}
 
 \emph{Let $X$ be an indecomposable complex in $\mathbf{C_{2}}(\mbox{proj} \,A)$ whose entries are $P_2\rightarrow P_1$. Consider $(\overleftarrow{X})_1: P_6\rightarrow P_2\oplus P_5 \rightarrow P_1$ and $(\overleftarrow{X})_2:P_4\oplus P_6\rightarrow P_2\oplus P_2 \rightarrow P_1$ two enlargements of $X$ in $\mathbf{C_{[0, 2]}}(\mbox{proj} \,A)$. It is not hard to see that  there is a morphism $f: (\overleftarrow{X})_1 \rightarrow (\overleftarrow{X})_2 $ such that $F(f) \simeq \left(
      \begin{array}{cc}
        1 & 0 \\
        0 & g \\
      \end{array}
    \right)$ where $g:(P_5\rightarrow 0)\longrightarrow (P_2\rightarrow 0) $ is a radical-irreducible morphism in $\mathbf{C_{2}}(\mbox{proj} \,A)$.}

\emph{On the other hand, the morphism $f$ is not irreducible in $\mathbf{C_{[0, 2]}}(\mbox{proj} \,A)$, since there is a non-trivial factorization $f=h_2h_1$ where $h_1:(\overleftarrow{X})_1\rightarrow Y $, $h_2:Y\rightarrow (\overleftarrow{X})_2$ and $Y$ is the indecomposable complex in $\mathbf{C_{[0, 2]}}(\mbox{proj} \,A)$ whose entries are $P_4\oplus P_6\rightarrow P_2\oplus P_5 \rightarrow P_1$. }

\end{ej}

The next result follows with the same techniques as the ones used in  \cite[Theorem 2.19]{BS}. For the convenience of the reader we state a proof.

\begin{lem}\label{projective} Let $A$ be an artin algebra. There are not  irreducible morphisms from $P$ to $P$ in $\emph{proj} \,A$, where $P$ is an indecomposable projective $A$-module if and only if there are not irreducible morphisms from $X$ to $X$ in $\emph{proj} \,A$ where $X$ is any projective $A$-module.
\end{lem}

\begin{proof} We only prove one implication  since the other is trivial. Assume that there is an irreducible morphism $f:X \rightarrow X$ where $X$ is any projective $A$-module. Moreover, assume that $X$ is not indecomposable, otherwise there is nothing to prove.

Suppose that $X=X_0 \oplus X_1$ with $X_0$ indecomposable in $\mbox{proj}\,A$ and $X_1$ a projective $A$-module. Then $f:X_0 \oplus X_1 \rightarrow X_0 \oplus X_1$.

$\it{Case} 1$. If $f$ is radical-irreduciblethen by Proposition \ref{BS3} we have that the morphism $g= f\sigma: X_0 \rightarrow X_0 \oplus X_1$ is irreducible, where $\sigma: X_0 \rightarrow X_0 \oplus X_1$ is a non-zero section. Hence $g$ is radical, since $X_0$ is indecomposable in $\mbox{proj}\,A$. Now, again by Proposition \ref{BS3} since $g: X_0 \rightarrow X_0 \oplus X_1$ is radical-irreducible then if we compose it with a non-zero retraction $ \pi: X_0 \oplus X_1  \rightarrow X_0$ we get that $\pi g \sigma:X_0 \rightarrow X_0$ is an irreducible morphism, which contradicts the hypothesis.

$\it{Case} 2$. Assume that $f$ is not radical. Then by Proposition \ref{BS1} we have that $f$ is isomorphic to a morphism of the form:
\begin{equation}\label{f1-f2} f=\left(  \begin{array}{cc}
       f_1 & 0 \\
        0 & f_2 \\
      \end{array}
    \right) : X_1 \oplus X_2 \rightarrow  Y_1 \oplus Y_2
    \end{equation}

\noindent with $f_1$ an isomorphism and $f_2$ a radical morphism. Since $f_1$ is an isomorphism we get $X_1 \simeq Y_1$ and thus $X_2\simeq Y_2$. Now, by Proposition \ref{BS2}, $f$ irreducible implies that $f_2$ is irreducible and  we get that $f_2$ is an radical-irreducible morphism from $X_2$ to $X_2$. Then we are again in  $\it{Case} 1$, proving the result.
\end{proof}

\vspace{.1in} In case that we consider $A$ a finite dimensional $k$-algebra over an algebraically closed field, it follows applying \cite[Lemma 4.2, Chapter I]{ASS} that there is an irreducible morphism from an indecomposable projective $A$-module $P$ to itself in $\mbox{proj} \,A$ if and only if $Q_A$ has a loop.
\vspace{.1in}

Next, we analyze how are the entries of an irreducible morphism between an extension of a complex $X \in \mathbf{C_{n}}(\mbox{proj} \,A)$ and an enlargement of $X$. For this purpose, first we prove the following result.  

\begin{teo}\label{enlarg-enlarg} Let $A$ be an artin algebra. Assume that there 
are not irreducible morphisms from $P$ to $P$ in $\emph{proj} \,A$, for any indecomposable projective $A$-module.
Let $f: X \rightarrow Y$ be an irreducible morphism in $\mathbf{C_{[0, n]}}(\emph{proj} \,A)$. Fix a direct sum decomposition of $X^0=Z_{0}\oplus Z_{1}$ and $Y^0=Z_{0} \oplus Z'_{1}$ such that where $Z_0$ is the direct sum of all common indecomposable direct summands of $X^0$ and $Y^0$ and  $Z_1$,$Z'_{1}$ have no non-zero common direct summands.   
 Then,
\begin{enumerate}
\item[(a)] $f^{0}$ is a section  if and only if $Z_{1}=0$. Moreover, if this is the case, $f$ is of the form \emph{(sec)}.
\item[(b)] $f^{0}$ is a retraction  if and only if $Z_{1}'=0$.
\item[(c)] $f^{0}$ is irreducible if and only if $Z_{1}\neq 0$ and $Z_{1}' \neq 0$. Moreover, if this is the case, $f$ is of the form \emph{(irr-sec)}.
\end{enumerate}
\end{teo}

 \begin{proof} (a) If  $f^{0}: Z_{0}  \oplus Z_{1} \rightarrow Z_0 \oplus Z_{1}'$ is a section then $Z_{1}$ is  a direct summand of $Z_{1}'$.
Since by hypothesis $Z_{1}$ and $Z_{1}'$ have no non-zero common direct summands then $Z_1=0$.

Assume that $Z_{1}=0$ and $Z_{0} \neq 0$. Then $f^{0}: Z_{0} \rightarrow Z_0 \oplus Z_{1}'$.  By Proposition \ref{GM}, we have three cases to analyze.

If $f^{0}$ is a section then we are done. Moreover, by Theorem \ref{GM}, since $f$ is irreducible, we have that $f^{i}$ is a section for all $i=0, \dots n$,
and we have that $f$ is of the form \mbox{(sec)}.

Suppose that $f^{0}$ is irreducible. Moreover, assume that $f^{0}$ is radical.
Let $\pi: Z_0 \oplus Z_{1}' \rightarrow Z'_0$ be any non-zero retraction. Then by Proposition \ref{BS3} we have
that $\pi f^{0}$ is an irreducible morphism from $Z_0$ to $Z_0$, which contradicts Lemma \ref{projective}. Therefore, $f^{0}$ is not radical.
By Proposition \ref{BS1}, we can assume that $f^{0}$ is isomorphic to a morphism of the form:

$$f^{0}=\left(
      \begin{array}{cc}
        1 & 0 \\
        0 & g \\
      \end{array}
    \right) : X_1 \oplus X_2 \rightarrow  X_1 \oplus Y_2$$

\noindent with $g$ a radical morphism,  where $Z_0= X_1 \oplus X_2$
and $Z_0 \oplus Z_{1}'= X_1 \oplus Y_2$. Moreover, we have that $X_2\neq 0$, otherwise $f^{0}$ is a section, contradicting our assumption.

On the other hand, $X_2$ is a direct summand of $Y_2$, that is, $Y_2=X_2\oplus Y_2'$. Then by Proposition \ref{BS2} $g: X_2 \rightarrow X_2 \oplus Y'_2$ is
 irreducible radical, because   $f^{0}$ is  irreducible. Therefore, $g$ is in the hypothesis of the previous case getting a contradiction to
 Lemma \ref{projective}.

Finally, assume that $f^{0}$ is a retraction. Then $Z_{1}'=0$. Therefore $f^{0}: Z_0 \rightarrow  Z_0$ is also a section.

%If for all $i$ we have that $f^{i}$ is a retraction then $W^{i}=0$ for all $i=0, \dots, n$ and
%we get that $f$ is an isomorphism a contradiction to the fact that $f$ is irreducible. Hence,  there is $j=0, \dots, n$ such
%that $f^{j}$ is irreducible in $\mbox{proj} \,A,$ getting a contradiction similar to the case below.

(b) The proof of Statement (b) is obtained by dual arguments of the ones used in the proof of Statement (a).

(c) It is an immediate consequence of the above statements and Proposition \ref{GM}. Moreover, by Proposition \ref{GM} we get that $f$ is
of the form \mbox{(irr-sec)}.
\end{proof}

Now, we shall give a method to decide the type of an irreducible morphism  any two indecomposable complexes in
$\mathbf{C_{n}}(\mbox{proj} \,A)$.
\vspace{.1in}

\noindent{\bf Method to decide the type of an irreducible morphism.}\label{met-enlarg-enlarg} Let $A$ be an artin algebra. Let $X,Y$ be  indecomposable complexes
in $\mathbf{C_{[0, n]}}(\mbox{proj} \,A)$. Assume that there are not irreducible morphisms from $P$ to $P$ in $\mbox{proj} \,A$, where $P$ is an indecomposable projective $A$-module.
 Let $f:X \to Y$ be an irreducible morphism in $\mathbf{C_{[0, n]}}(\mbox{proj} \,A)$. We fix a decomposition of each entry of the given complexes as in our previous theorem follows:

\noindent $X:X^0\oplus Z^ 0 \rightarrow X^{1}\oplus Z^{1} \rightarrow X^{2}\oplus Z^{2}\rightarrow \dots \rightarrow X^{n}\oplus Z^{n}$ and
$Y: X^0 \oplus W^0  \rightarrow X^{1} \oplus W^{1} \rightarrow X^{2} \oplus W^{2}\rightarrow \dots \rightarrow X^{n}\oplus W^{n}$
where $Z^ {i}$ and $W^{i}$ have no non-zero common direct summands.

Then we proceed analyzing each entry and we conclude as follows:
\begin{enumerate}
\item[(1)] If $Z^{0}=0$  then $f^{0}$ is a section and therefore $f$ is of the form \mbox{(sec)}.

\item[(2)] Assume $Z^0 \neq 0.$ If $W_{0}=0$ then $f^{0}$ is a retraction (and not a section). Then we have two possibilities, either $f$ is of type \mbox{(ret)} or $f$ is of the form \mbox{(ret-irr-sec)}.  If for all $i$, $i=1, \dots, n$, we have that $W^{i}=0$ then $f^{i}$ is a retraction for all $i$ and we get that $f$ is of the form \mbox{(ret)}. In any other case, we have that $f$ is of the form 
 \mbox{(ret-irr-sec)}.
 \item[(3)] If $Z_{0}\neq 0$ and $W_{0} \neq 0$  then $f^{0}$ is irreducible and therefore $f$ is of the form \mbox{(irr-sec)}.
\end{enumerate}

\vspace{.1in}

As an application of our previous method we prove the following generalization of \cite[Proposition 4.1]{CGP6}

\begin{prop}\label{Tipo}  Let $A$ be an artin algebra. Assume that there 
are not irreducible morphisms from $P$ to $P$ in $\emph{proj} \,A$, for any indecomposable projective $A$-module. Let $f:X \rightarrow Y$ and $g:X \rightarrow Y$ be irreducible morphisms between indecomposable complexes in $\mathbf{C_n}(\rm { proj}\,A)$.
Then the following statements hold.
\begin{enumerate}
\item [(a)] $f$ is of the form $(sec)$ if and only if $g$ is of the form $(sec)$.
\item [(b)] $f$ is of the form $(ret)$ if and only if $g$ is of the form $(ret)$.
\item [(c)] $f$ is of the form $(ret-irr-sec)$ with $f^i$ irreducible in $\rm{ proj}\, A$ if and only if $g$ is of the form $(ret-irr-sec)$ with $g^i$ irreducible in $\rm{ proj}\, A$.
\end{enumerate}
\end{prop}

As an immediate consequence of theorem \ref{enlarg-enlarg}, we have the following corollary.

\begin{coro}\label{ext-enlarg} Let $A$ be an artin algebra. Let $X$ be an indecomposable complex in $\mathbf{C_{n}}(\emph{proj} \,A)$.
Assume that there are not irreducible morphisms from $P$ to $P$ in $\emph{proj} \,A$, where $P$ is an indecomposable projective $A$-module.
Let $f: {(\overleftarrow{X})_1} \rightarrow {^{\epsilon}X}$ be an irreducible morphism in $\mathbf{C_{[0, n]}}(\emph{proj} \,A)$,
where $(\overleftarrow{X})_1$ and $^{\epsilon}X$ are an enlargement which is not an extension and an extension of $X$, respectively. Assume that $(\overleftarrow{X})_1$ and $^{\epsilon}X$ are as follows:

\noindent $(\overleftarrow{X})_1: X_0 \oplus Z_{0} \rightarrow X^{1} \oplus Z^{1} \rightarrow X^{2} \oplus Z^{2}\rightarrow \dots \rightarrow X^{n}\oplus Z^{n}$ and $^{\epsilon}X : X_0 \oplus Z'_1 \rightarrow X^{1} \rightarrow X^{2} \rightarrow \dots \rightarrow X^{n}$.  Then we have that $Z'_1=0$ and $f$ is of the form \emph{(ret)}.
\end{coro}

\begin{proof} Since $(\overleftarrow{X})_1$ is an enlargement and not an extension then $Z^1 \neq 0$, otherwise $X$ decomposes.  Therefore by (2) we have that $f^{1}$ is a retraction and not a section. Moreover, in this case $f^{0}$ is also a  retraction then $Z'_1=0$. We observe that for all $i$,  $f^{i}$ is not irreducible. Then for all $i$,  $f^{i}$ is a retraction, and there have to be $j=1, \dots, n$ where $Z^{j} \neq 0$, since $f$ is irreducible. Therefore $f$ is of the form \mbox{(ret)}.
\end{proof}

Next, we illustrate in an example how to apply the above method.

\begin{ej}
\emph{ Let $A$ be the algebra given by the quiver:}

{\begin{displaymath}
\xymatrix @R=0.1cm  @C=1cm{1 \ar[rd]^{\alpha_1}& & & & & &  8  \\ & 3\ar[r]^{\alpha_2} & 4\ar[r]^{\alpha_3}& 5 \ar[r]^{\alpha_4}& 6 \ar[r]^{\alpha_5} & 7  \ar[rd]_{\beta_2}\ar[ru]^{\alpha_6}&  \\
2 \ar[ru]_{\beta_1}& & & & & &   9}
\end{displaymath}}

\noindent \emph{with the relations $\alpha_2\alpha_1= \alpha_2 \beta_1 = \alpha_4 \alpha_3 = \beta_2 \alpha_5 =\alpha_6 \alpha_5=0$.}

\emph{Since $A$ is a piecewise hereditary algebra there are not irreducible morphisms from $P$ to $P$ in $\mbox{proj} \,A$, where $P$ is an indecomposable projective
$A$-module. By constructing the preprojective component of the Auslander-Reiten quiver of $\mathbf{C_{[0, 2]}}(\mbox{proj} \,A)$, we can see that
there is an irreducible morphism $f$  between two enlargements of an indecomposable
complex in $X \in \mathbf{C_{[0, 2]}}(\mbox{proj} \,A)$ as follows:}
\vspace{.1in}

\begin{displaymath}
\xymatrix { P_9 \oplus P_8 \ar[r]\ar[d]^{f^{0}}& P_7 \oplus P_7  \ar[r]\ar[d]^{f^{1}}
&  P_5 \oplus P_6  \ar[r]\ar[d]^{f^{2}} & P_3  \ar[d]^{f^{3}}\\
 P_9 \oplus P_8 \ar[r]& P_7 \oplus P_7  \ar[r]
&  P_5 \oplus P_5  \ar[r]  & P_4 \oplus P_3}
\end{displaymath}

\noindent \emph{where $X$ is the indecomposable complex whose entries are  $P_7 \rightarrow P_5 \rightarrow P_3$. }

\emph{ Reconsidering the entries of the morphism $f$ as follows:}

\begin{displaymath}
\xymatrix { (P_9 \oplus P_8)\oplus 0 \ar[r]\ar[d]^{f^{0}}& (P_7 \oplus P_7) \oplus 0  \ar[r]\ar[d]^{f^{1}}
&  (P_5) \oplus (P_6)  \ar[r]\ar[d]^{f^{2}} & (P_3) \oplus 0 \ar[d]^{f^{3}}\\
(P_9 \oplus P_8)\oplus 0 \ar[r]& (P_7 \oplus P_7)\oplus 0  \ar[r]
&  (P_5) \oplus (P_5)  \ar[r]  & (P_3) \oplus (P_4)}
\end{displaymath}

\noindent \emph{and applying the above method, %\ref{met-enlarg-enlarg},
we can observe that $f^{2}$ is an irreducible morphism in $\mbox{proj} \,A$.
Therefore $f$ is of the form (\mbox{ret-irr-sec}).}
\end{ej}

%\begin{teo}\label{ext-enlarg-general} Let $A$ be an artin algebra. Let $X$ be an indecomposable complex in $\mathbf{C_{n}}(\emph{proj} \,A)$.
%Assume that there are not irreducible morphisms from $P$ to $P$, where $P$ is an indecomposable projective $A$-module.
%Let $f: {\epsilon}_X \rightarrow \overleftarrow{X}$ be an irreducible morphism in $\mathbf{C_{n+1}}(\emph{proj} \,A)$,
%where the indecomposable complexes ${\epsilon}_X$ and  $\overleftarrow{X}$ are as follows:
%
%\noindent ${\epsilon}_X: X_0 \rightarrow X^{1} \rightarrow X^{2}\rightarrow \dots \rightarrow X^{n}$ and
%$\overleftarrow{X}: Z^{0} \rightarrow X^{1} \oplus Z^{1} \rightarrow X^{2} \oplus Z^{2}\rightarrow \dots \rightarrow X^{n}\oplus Z^{n}$.
%
%Then $f$ is of the form  VER.
%\end{teo}

\vspace{.1in}

The authors have no financial or proprietary interest in the material discussed in this article.

\end{document}